\documentclass[sn-mathphys-num]{sn-jnl}

\usepackage{graphicx}%
\usepackage{multirow}%
\usepackage{amsmath,amssymb,amsfonts}%
\usepackage{amsthm}%
\usepackage{mathrsfs}%
\usepackage[title]{appendix}%
\usepackage{xcolor}%
\usepackage{textcomp}%
\usepackage{manyfoot}%
\usepackage{booktabs}%
\usepackage{algorithm}%
\usepackage{algorithmicx}%
\usepackage{algpseudocode}%
\usepackage{listings}%

\newtheorem{theorem}{Theorem}
\newtheorem{remark}{Remark}%

\newtheorem{lemma}{Lemma}%

\begin{document}

\title[Hybrid CGME and TCGME algorithms]
{Hybrid CGME and TCGME algorithms for
large-scale general-form regularization}


\author*[1]{\fnm{Yanfei} \sur{Yang}}\email{yang@zafu.edu.cn}

%

\affil*[1]{\orgdiv{School of mathematics and computation science}, \orgname{Zhejiang A\&F University}, \city{Hangzhou}, \postcode{300113}, \country{China}}

%


\abstract{
Two new hybrid algorithms are proposed for large-scale linear
discrete ill-posed problems in general-form regularization.
They are both based on Krylov subspace inner-outer iterative algorithms.
At each iteration, they need to solve a linear least squares problem, which is the
inner least squares problem.
It is proved that inner linear least squares problems,
solved by LSQR, become better conditioned as $k$ increases,
so LSQR converges faster.
We also prove how to choose the stopping tolerance
for LSQR to guarantee that the computed
and exact best regularized solutions have the same accuracy.
Numerical experiments are provided
to demonstrate the effectiveness and efficiency of our new hybrid algorithms, along with comparisons to the existing algorithm.}

\keywords{Linear discrete ill-posed,
General-form regularization, hybrid algorithms, 
CGME, the truncated CGME algorithm}


\pacs[MSC Classification]{65F22, 65F10, 65J20, 65F35, 65F50}

\maketitle

\section{Introduction}\label{sec1}

Consider the large-scale linear discrete ill-posed problem of the form
\begin{equation}\label{eq1}
\min_{x\in \mathbb{R}^{n}}\|Ax - b\| \  {\rm or}
\ Ax = b, \ A \in \mathbb{R}^{m \times n}, \ b \in \mathbb{R}^{m},
\end{equation}
where the norm $\|\cdot\|$ is the 2-norm of a vector or matrix,
and $A$ is ill-conditioned with its singular values
generally decaying to zero without
a noticeable gap between consecutive ones,  
and the right-hand side
$b = b_{true}+e$ is assumed to be contaminated
by a Gaussian white noise $e$, where $b_{true}$
is the noise-free right-hand side and $\|e\| < \|b_{true}\|$.
Without loss of generality, assume that
$Ax_{true} = b_{true}$.
Discrete linear ill-posed problems of the form (\ref{eq1}) are
derived from the discretization of linear ill-posed problems,
such as the first kind Fredholm integral equation,
and arise in various scientific research and applications areas,
including biomedical sciences, geoscience, mining engineering and
astronomy; see, e.g.,
\cite{aster,berisha2012restore,engl96,epstein2007,haber2014,hansen10,
kirsch11,miller,nat01,vogel02}.

Because of the presence of the noise $e$ and the high ill-conditioning
of $A$, the naive solution
$x_{naive}=A^{\dag}b$ to (\ref{eq1}) generally
is meaningless since it
bears no relation to the true solution
$x_{true}=A^{\dag}b_{true}$,
where $\dag$ denotes the Moore-Penrose inverse of a matrix.
To calculate a meaningful solution,
it is necessary to employ regularization
to overcome the inherent instability of ill-posed
problems. The basic idea of regularization is that the underlying problem
(\ref{eq1}) is replaced with a modified problem which is
relatively stable and can
obtain a regularized solution to approximate the true solution.
There are various regularization techniques taking many forms;
see, e.g., \cite{chung,hansen98,hansen10}.

A simple and popular regularization
is by iterative regularization.
In this setting,
an iterative method is applied directly to
\begin{equation*}
\min_{x\in \mathbb{R}^{n}}\|Ax - b\|
\end{equation*}
and regularization is obtained by
terminating early.
It is well known that
iterative algorithms exhibit semi-convergence
on ill-posed problems, with
errors
initially decreasing but at some point beginning to
increase since the small singular values of $A$ start to amplify noise~\cite{hansen10, hansen98, jia20203}.
Therefore, for iterative regularization,
the number of iterations plays the role of regularization parameter and
a vital and nontrivial
task is to select a good stopping iteration,
which means determining a good regularization parameter.
Unfortunately, some iterative algorithms
only have partial regularization, which means
only an iterative algorithm can not get the best regularized solution; see, e.g., \cite{jia20203, jia2020low,jia2020lsqr} for more details.

Another popular way is the {hybrid regularization} method.
Since O'Leary and Simmons~\cite{oleary1981} introduce
the hybrid method using the Golub-Kahan bidiagonalization
process and truncated singular value decomposition (SVD) for large-scale standard-form regularization,
various hybrid methods based on the
Krylov subspace method have been introduced
including more general-form regularization
terms and constraints; see, e.g.,  \cite{ chung19,chung,chung15,chung17,
gazzola2014,gazzola2014(2),gazzola2015,gazzola2019,
lampe12, bazan2014}.
The hybrid regularization method
generally calculates an approximation of $x_{true}$
by solving the following problem
\begin{equation}\label{gen1}
\min_{x\in\mathbb{S}}\|Lx\| \ {\rm subject \ to } \
\mathbb{S}=\left\{x|\|Ax-b\|\leq\tau \|e\|\right\}
\end{equation}
or
\begin{equation}\label{gen2}
\min_{x\in\mathbb{R}^n}\left\{\|Ax-b\|^2+\lambda^2\|Lx\|^2\right\},
\end{equation}
where $\tau\approx 1$,
$L\in \mathbb{R}^{p\times n}$
is a regularization matrix and
usually a discrete approximation of some derivative operators,
and $\lambda>0$ is the regularization parameter
that controls the amount of regularization to balance
the fitting term and the regularization term; see, e.g.,
\cite{hansen98,hansen10,tikhonov63}.
Problem (\ref{gen1}) is the discrepancy principle-based general-form regularization which is equivalent to Problem (\ref{gen2}) which is the general-form Tikhonov regularization;
see~\cite[pp.63-4, 172, 181-2]{hansen10}
and~\cite[pp.11, 85, 105, 179]{hansen98}
for the equivalence of the above two formulations.
The solution to (\ref{gen2}) is unique for a given $\lambda>0$ when
$N(A)\cap N(L)={0}$
where $N(\cdot)$ denotes the null space of a matrix.
When $L = I_n$, with $I_n$ being the $n \times n$ identity matrix,
Problem (\ref{gen1}) and Problem (\ref{gen2})
are said to be in standard form.


It is worth mentioning that regularization in other norms than the 2-norm is also significant,
and some ill-posed problems may have an underlying
mathematical model that is not linear.
Consider the general optimizations of the form
\begin{equation*}\label{gen_1}
\min_{x\in \mathbb{S}}\mathcal{R}(Lx)  \ {\rm subject \ to } \
\mathbb{S }=\left\{x|\mathcal{J}(Ax-b)= \min \right\}
\end{equation*}
or
\begin{equation*}\label{gen_2}
\min_{x\in\mathbb{R}^n}\{\mathcal{J}(Ax-b)+\lambda^2\mathcal{R}(Lx)\},
\end{equation*}
where $\mathcal{J}$ is a fit-to-data term and $\mathcal{R}$ is a regularization term.
It is known, for instance, that solving
\begin{equation}\label{gen_3}
\min_{x\in\mathbb{S}}\|Lx\|_p^p \  {\rm subject \ to } \
\mathbb{S}=\{x|\|Ax-b\|_q^q =\min \}
\end{equation}
or
\begin{equation}\label{gen_4}
\min_{x\in\mathbb{R}^n}\{\|Ax-b\|_q^q+\lambda^2\|Lx\|_p^p\},
\end{equation}
where $1\leq q<2$ and $1\leq p<2$.
For these cases, one should take advantage of more sophisticated methods to solve nonlinear optimization problems;
however, many of these methods require solving a subproblem
with an approximate linear model;
see, e.g., \cite{chung19,gazzola2014} and \cite[pp.~120-121]{hansen98} for more details.

In this paper, we focus our discussion on the linear problem, i.e., Problem (\ref{gen1}).
The reason we focus our study on Problem (\ref{gen1}) instead of Problem (\ref{gen2}) is that
 when we solve (\ref{gen1}), the number of iterations plays the regularization parameter,
which means that
we have not to determine the optimal regularization parameter prior to the solution at each iteration.
The basic idea of our new hybrid method is first to solve the underlying problem (\ref{eq1}) using a Krylov solver
and then to apply a general form regularization term to the projected problems generated by the Krylov solver.
Our method means that only the constraint of Problem (\ref{gen1}) is projected
onto a consequence of nested Krylov subspaces, while the
regularization part remains unchanged.
This approach, on the one hand, can make full use of the advantages of the Krylov subspace method to reduce the large-scale problem to a small or medium problem.
On the other hand, it can also retain as much information as possible.
Therefore, our method expects to capture
as much information as possible about the dominant generalized singular value decomposition (GSVD)
of the matrix pair $\{A, L\}$ components.

Before proceeding, it should be emphasized that regardless of the method used, the most important basic ingredient for solving (\ref{gen1}) and (\ref{gen2}) successfully is that the best regularized solution must capture the dominant GSVD components of the matrix pair $\{A, L\}$ while suppressing those corresponding to small generalized singular values; see, e.g., \cite{hansen98, hansen10, kilmer07, jia2020}.

Kilmer \textit{et al.} \cite{kilmer07}
adapt the joint bidiagonalization (JBD) process
and develop a JBD process that successively reduces the matrix pair $\{A, L\}$ to lower and upper bidiagonal forms, respectively.
Based on this process, they propose an iterative algorithm based on projections to solve (\ref{gen2}).
It is argued \cite{kilmer07} that the underlying
solution subspaces are legitimate since they appear to be more directly related to the generalized right singular vectors of the pair of matrix $\{A, L\}$ solution subspace.
However, in the method, one needs to determine an optimal or suitable
regularization parameter for each projected general-form Tikhonov regularization problem, and, at the same time, one has to
judge when the projection subspace is large enough. Unfortunately, the meaning of 'sufficiently large' or 'large enough' has been vague, and no deterministic determination has been given without ambiguity up to now.
Jia and Yang \cite{jia2020} propose a new JBD process-based iterative algorithm, called the JBDQR algorithm, for solving (\ref{gen1}) instead of (\ref{gen2}). The JBDQR algorithm uses the same JBD process as the one due to Kilmer \textit{et al.} \cite{kilmer07}. Therefore, they have the same underlying solution subspace that is argued to be legitimate.
Although they are based on the same JBD process and have the same underlying solution subspace,
the mechanism and features of the two algorithms are different and it is argued in \cite{jia2020}.

The JBD process is an inner-outer iterative process.
At each outer iteration, the JBD process needs to compute the solution of a large-scale linear least squares problem with the coefficient matrix $(A^T, L^T )^T$ that is larger than the problem itself and is supposed to be solved iteratively, called inner iteration. Fortunately, $(A^T, L^T )^T$ is generally well conditioned, as $L$ is typically in applications \cite{hansen98,hansen10}.
In these cases, the LSQR algorithm \cite{lsqr1982} can solve the least squares problems mentioned efficiently.
 Although the underlying
solution subspaces generated by the process are legitimate,
the overhead of the methods based on the process may be extremely expensive.
Unfortunately, methods based on the
JBD process cannot avoid solving the large-scale inner least squares problem at every iteration. Finally, methods based on the process need to solve a large-scale least squares problem with the coefficient matrix $(A^T , L^T )^T$ to form a regularized solution. Therefore, the method based on the JBD process may
come at considerable computational costs.

Novati and Russo \cite{novati2014gcv} propose a Arnoldi-Tikhonov (AT)
method for solving (\ref{gen2}). First, the method uses the Arnoldi process to
reduce the underlying matrix $A$ to a sequence of small upper Hessenberg matrix. Then
the regularization matrix $L$ is multiplied by the column orthogonal matrix which is generated by the Arnoldi process applied to $A$ with starting vector $b$.
In other words, the AT method first projects $A$ onto a sequence of Krylov subspaces, while
projects $L$ onto the same Krylov subspaces that have nothing to do with $L$ because they are generated only by $A$ and $b$.
Gazzola and Novati \cite{gazzola2013multi} propose an iterative algorithm based on the AT method and the discrepancy principle \cite{hansen10,hansen98} for multiparameter Tikhonov regularization problems.
They also present a generalized AT (GAT) method in \cite{gazzola2014(2)} for solving (\ref{gen2}).
The difference between the GAT method and
the AT method is the parameter selection method.
Gazzola and Nagy \cite{gazzola2014} applying the flexible Krylov subspaces and the restarting Arnoldi algorithm to the AT method obtain two new algorithms for solving (\ref{gen2}).
Novati and Russo \cite{novati2014adaptive} propose an
adaptive AT algorithm for image restoration.

These studies about AT and GAT methods are all based on the Arnoldi process.
The process requires that the underlying matrix $A$ must be a square matrix, and
the process takes advantage of only $A$ but does not use its transpose, so the
information on $A^T$ is lacking.
In fact, for $L=I_n$, Hansen \cite[p.126]{hansen10} points out that
the success of Arnoldi process-based methods is highly problem dependent,
because they mix the SVD components in each iteration,
and they can be successful when the mixing of the SVD components is weak, e.g.,
$A$ is (nearly) symmetric.
Therefore, these algorithms are not theoretically guaranteed to capture the dominant SVD components of $A$.
Furthermore, the regularization matrix $L$ is blindly projected onto the Krylove subspaces generated only by $A$ and $b$, therefore, the algorithms may lose important information on the regularization term.
Consequently, there is no guarantee that these algorithms can capture the dominant GSVD components of the matrix pair $\{A,L\}$.
This may cause the algorithms to fail, or even if the algorithms work, the accuracy of the regularized solution obtained by the algorithms may not be so high because, as we mentioned above, the key to solving (\ref{gen1}) and (\ref{gen2}) successfully
is that the regularized solution must capture all the needed dominant GSVD components of the matrix pair $\{A, L\}$; see, e.g., \cite{hansen98,hansen10,jia2020, kilmer07}.
In addition, these algorithms are all for solving
Problem (\ref{gen2}). These algorithms need to determine the
optimal regularization parameter prior to the solution.
In fact, it is hard to prove the optimal regularization
parameter of the projected problem is the one of Problem (\ref{gen1}).

Gazzola \textit{et al.} \cite{gazzola2015} propose an iterative algorithm based on
the nonsymmeteric Lanczos process solving (\ref{gen2}).
Viloche Bazn \textit{et al.} \cite{viloche2014extension} extend the
GKB-FP algorithm proposed by Bazán and
Borges \cite{bazan2014} to general-form Tikhonov regularization.
They carry the GKB-FP algorithm over to standard-form transformation
and free-of-standard-form transformation, respectively.

All of these studies are based on the Golub-Kahan bidiagonazation
process.
These algorithms, except for the extension of GKB-FP for standard-form transformation,
first project the underlying matrix $A$
onto a sequence of Krylov subspaces and
then simultaneously project the regularization matrix $L$
onto the same Krylov subspaces which only generated by $A$ and $b$.
Finally, QR factorization or GSVD and parameter choice methods are
used to solve the projected problems.
Although Gloub-Kahan process based algorithms, such as
LSQR and CGME, has partial or full regularization \cite{jia2020lsqr, jia20203,
 jia2020low} which means they can capture some or all of the dominant SVD components of $A$,
blindly projecting the regularization matrix $L$ onto the Krylov subspaces, just as we mentioned above,
may lost important information about the regularization term. As a consequence,
these algorithms can not be theoretically guaranteed to capture the dominant GSVD components of $\{A,L\}$, which is crucial to solve (\ref{gen1}) and (\ref{gen2}) successfully.

In this paper,
we propose two new hybrid algorithms for solving (\ref{gen1}) instead of (\ref{gen2}).
We first  exploit a Krylov solver to solve the underlying problem (\ref{eq1}).
Then a general-form regularization term is applied to the projected problems
which are generated by the Krylov solver.
Finally, we solve the projected problem with the general-form regularization term.
They are both inner-outer iterative algorithms.
At every iteration, we need to solve a linear least squares problem,
called inner iteration.
The resulting algorithms are called
hybrid CGME and hybrid truncation CGME (TCGME),
because the Krylov solvers are the CGME
algorithm and the TCGME algorithm presented by Jia in \cite{jia20203}, respectively.
They are abbreviated as hyb-CGME
and hyb-TCGME, respectively.

One of the benefits of our method is that it can capture the dominant SVD components of $A$
and simultaneously it also includes all the information of $L$,
therefore, it is can be expected that our
method can capture the dominant GSVD of the matrix pair $\{A, L\} $ components information as much as possible, which is the key to solve (\ref{gen1}) and (\ref{gen2}) successfully.
Another benefit is that we no longer need to determine an optimal regularization parameter for each projected problem with general-form regularization prior to the solution, which itself, unlike the optimal regularization parameter of the original problem (\ref{gen2}), may be hard to define because the projected problem may fail to satisfy the discrete Picard condition, so that regularized solutions may behave irregularly.
The other benefit is that we avoid solving the
large-scale linear least squares problems with
the coefficient matrix $(A^T,L^T)^T$ in JBDQR to
obtain the regularized solution,
and the size of inner linear least squares problems is much less than the one of JBDQR.
Therefore, it is can be expected that our method is substantially cheaper than JBDQR.

We will provide strong theoretical support for our algorithms by establishing a number of results.
we will prove that the inner least squares problems become better conditioned as $k$ increases.
Then LSQR chosen to solve the inner least squares problems converges faster as $k$ increases.
In principle, the inner linear least squares problems are supposed to be solved accurately,
but in practice, they are solved by an iterative algorithm
since they are large-scale problems.
We will prove how to choose the stopping
tolerance for the iterative algorithm
to solve the inner least squares problems,
To guarantee that the computed regularized
solution has the same accuracy as the accurate regularized solution.


The rest of this paper is organized as follows.
In Section 2, we briefly review
CGME and TCGME.
We propose hyb-CGME and hyb-TCGME algorithms and
make an analysis on the
conditioning of inner least squares problems in Section 3.
In Section 4, we make a theoretical analysis on the stopping tolerance for LSQR.
Numerical experiments are presented in Section 5.
Finally, we conclude the paper in Section 6.

\section{ CGME and truncated CGME algorithms}\label{sec2}

In this section we provide some necessary background. We describe
the CGME and truncated CGME (TCGME) algorithms which are
based on the Golub-Kahan bidiagonalization process.
Given the initial vectors $\beta_1p_1=b$ and $\beta_1q_0=0$,
for $i=1, 2, 3, \cdots$, the Golub-Kahan
iterative bidiagonalization computes
\begin{eqnarray}
\alpha_{i}q_{i}&=&A^Tp_{i}-\beta_{i}q_{i-1}, \label{bid1}\\
\beta_{i+1}p_{i+1}&=&Aq_i-\alpha_ip_{i},\label{bid2}
\end{eqnarray}
where $\beta_{i+1}\geq 0$ and $\alpha_{i+1}\geq 0$
are normalization constants chosen so that $\|q_i\|=\|p_{i+1}\|=1$.
In particular, $\beta_1=\|b\|$.

With the definitions
\begin{equation}\label{lsqr1}
Q_k=(q_1, \ldots,q_k),\ P_k = (p_1,\ldots,p_k),
\end{equation}
and
\begin{equation}\label{bk}
B_k = \left(
\begin{array}{cccc}
\alpha_1 &  &  &  \\
\beta_2 & \alpha_2 &  &  \\
& \ddots &\ddots  &  \\
&  & \beta_k & \alpha_k \\
\end{array}
\right)\in \mathbb{R}^{k\times k},
\ B_{k+}=\left(
       \begin{array}{c}
         B_k \\
         \beta_{k+1}e_k^T \\
       \end{array}
     \right)\in \mathbb{R}^{{(k+1)}\times k},
\end{equation}
where $e_k$ is the $k$-th vector of the standard Euclidean basis vector of
$\mathbb{R}^{k}$,
the recurrence relations (\ref{bid1}) and (\ref{bid2})
can be written of the form
\begin{eqnarray*}
  A^{T}P_{k}&=&Q_{k}B_k^T,\label{gol}\\
   AQ_k&=&P_{k+1}B_{k+}.
\end{eqnarray*}
These matrices can be computed by the following $k$-step Golub-Kahan bidiagonalization process.

\begin{algorithm}
\caption{\ $k$-step Golub-Kahan bidiagonalization process.}
\label{alg1}
\begin{algorithmic}[1]
\State Take $p_1 = b/\|b\|\in \mathbb{R}^{m}$, and define $\beta_1q_0 = 0$.
\For{$j = 1, 2, \cdots, k$ }
\State $r = A^Tp_j - \beta_jq_{j-1}$
\State $\alpha_j = \|r\|$; ~$q_j = r/\alpha_j$
\State  $s = Aq_j - \alpha_jp_j$
\State  $\beta_{j+1} = \|s\|$; ~$ p_{j+1} = s/\beta_{j+1}$.
\EndFor
\end{algorithmic}
\end{algorithm}

It is well known that CGME \cite{bjorck,hanke,hnety,jia20203}
is the CG method implicitly applied to
\begin{equation*}
\min\|AA^Ty-b\| \ {\rm or} \ AA^Ty=b \ {\rm and} \ x=A^Ty,
\end{equation*}
and it solves the problems
\begin{equation*}
\|x_{naive}-x_k^{cgme}\|=\min_{x\in\mathcal{V}_k}\|x_{naive}-x\|
\end{equation*}
for the iterate $x_k^{cgme}$, where $\mathcal{V}_k=\mathcal{K}_k(A^TA,A^Tb)$
is the $k$
dimensional Krylove subspace generated by the
Golub-Kahan bidiagonalization process.

From the $k$-step Golub-Kahan
bidiagonalization process, it follows
\begin{equation}\label{pro_2}
{B}_{k}=P_{k}^TAQ_k.
\end{equation}
Noting $\|b\|e_1=P_k^Tb$ and (\ref{pro_2}), we have
\begin{equation}\label{xk}
x_k^{cgme}=Q_k{B}_k^{-1}P_k^Tb.
\end{equation}
Therefore, CGME solves a sequence of problems
\begin{equation*}
\min \|P_k{B}_kQ_k^Tx-b\|
\end{equation*}
for $x_k^{cgme}$ starting with $k=1$,
where the projection $P_k{B}_kQ_k^T$, a rank-$k$ approximation to $A$,
substitutes $A$ in the underlying ill-posed problem (\ref{eq1}).

TCGME \cite{jia20203} solves a sequence of problems
\begin{equation}\label{tc}
\min \|P_{k+1}{C}_kQ_{k+1}^Tx-b\|
\end{equation}
 for $x_{k}^{tcgme}$ starting with $k=1$, where ${C}_k$ is the best rank-$k$ approximation for ${B}_{k+1}$ defined in (\ref{bk}).
Obviously, the solution to (\ref{tc}) is
\begin{equation}\label{txk}
x_{k}^{tcgme}=Q_{k+1}{C}_k^{-1}P_{k+1}^Tb.
\end{equation}

About the accuracy of the rank-$k$ approximations
$P_{k}{B}_kQ_k^T$  and $P_{k+1}{C}_kQ_{k+1}^T$ to $A$,
Jia has established the following results (cf. \cite[Theorem 1 and Theorem 5]{jia20203}).

\begin{theorem}\label{jia20201}
For the rank-$k$ approximations
$P_k{B}_kQ_k^T$ and $P_{k+1}{C}_kQ_{k+1}^T$ to $A$,
$k=1,2,\dots,n-1,$ we have
\begin{eqnarray}
\gamma_k^{lsqr}<\gamma_k^{cgme}<\gamma_{k-1}^{lsqr},\label{cgme1}\\
\gamma_{k+1}^{cgme}<\gamma_k^{cgme},\label{cgme}\\
\gamma_{k}^{tcgme}\leq{\theta}_{k+1}^{k+1}+\gamma_{k+1}^{cgme}\label{tcgme},
\end{eqnarray}
where ${\theta}_{k+1}^{k+1}$ is the smallest singular value of ${B}_{k+}$,
$\gamma_k^{cgme}=\|A-P_k{B}_kQ_k^T\|, \gamma_k^{tcgme}=\|A-P_{k+1}{C}_kQ_{k+1}^T\|,
\gamma_k^{lsqr}=\|A-P_{k+1}B_{k+}Q_k^T\|$, and $\gamma_0^{lsqr}=\|A\|$.
\end{theorem}

Jia in \cite{jia20203} points that inequalities (\ref{cgme1}) imply that
$P_{k}{B}_kQ_k^T$ is definitely a less accurate rank-$k$
approximation to $A$ than $P_{k+1}B_{k+}Q_k^T$ and
there is no guarantee that $P_{k}{B}_kQ_k^T$
is a near best rank-$k$ approximation
to $A$ even for severely and moderately ill-posed problems.
From (\ref{cgme}) and (\ref{tcgme}), it follows that
$$
\gamma_k^{tcgme}<{\theta}_{k+1}^{k+1}+\gamma_{k}^{cgme}.
$$
The above relation implies that
when ${\theta}_{k+1}^{k+1}$ is very small, the rank-$k$ approximation $P_{k+1}C_kQ^T_{k+1}$
to $A$ is more accurate than
the rank-$k$ approximation $P_{k}B_{k}Q^T_k$.

As there
is a lack of theory to guarantee that
$P_kB_kQ^T_k$ and $P_{k+1}C_{k}Q^T_{k+1}$
are the best rank-$k$
approximations to $A$, in fact,
it is not even theoretically guaranteed that
$P_kB_kQ^T_k$ and $P_{k+1}C_{k}Q^T_{k+1}$
are the near best rank-$k$
approximations to $A$,
which means the regularized solution obtained by Krylov solvers, CGME or TCGME, alone possibly
only has partial regularization \cite{jia20203}.
In other words, the regularized solution obtained
by CGME or TCGME alone at the semi-convergence point may not be the
best one.
We now recall the semi-convergence phenomenon:
At the beginning of the iterative process, the iterates converge to $x_{true}$;
afterwards, the noise $e$ starts to deteriorate
the iterates so that they start to diverge from $x_{true}$ and instead converge to $x_{naive}$,
which is a well-known basic property of Krylov solvers, see, e.g., \cite{hansen98, hansen10, gazzola2015}.
Semi-convergence is due to the fact that
the projected problem starts to inherit the ill-conditioning of
(\ref{eq1}) after a certain number of steps.
That is to say, before the best regularized solution is obtained,
the projected problems generated by Krylov solvers, CGME and TCGME,
start to inherit the ill-condution of (\ref{eq1}).
To get the best regularized solution
which is as accurate as the regularized solution with full regularization,
we consider to employ the general-form regularization term to the
projected problems,
which is the basic idea of the hybrid method.

\section{New hybrid algorithms}\label{sec3}

 According to the above analysis, 
To obtain the best possible regularized solution,
we shall propose two projection based hybrid algorithms which solve (\ref{gen1}).
Our method would be as easy to implement as the JBDQR algorithm 
and much cheaper than it because at each iteration it solves 
an inner linear least squares problem, 
the size of which is significantly smaller than that of JBDQR,
and our method avoids solving a large-scale
linear least squares problem to obtain the regularized solution;
see \cite{jia2020} for more details.
We will establish a number of theoretical results and get insight into the effectiveness of the method.
Particularly, we will prove how the conditioning of the inner least squares problems changes as $k$ increases and prove how to choose the stopping tolerance of LSQR which is chosen to solve the inner least squares problems.

Now, we propose the new hybrid algorithms.
First, the Krylov solvers, CGME and TCGME, are used to solve the underlying problem, respectively.
This means the underlying problem is projected onto
\begin{equation}\label{pro1}
\left\|P_{k}{B}_kQ_k^Tx-b\right\| = \min
\end{equation}
 or
\begin{equation}\label{pro2}
\left\|P_{k+1}{C}_kQ_{k+1}^Tx-b\right\| = \min,
\end{equation}
where $k=1,2,\cdots,n.$ With the iteration,
before the semi-convergence point,
the projections $P_{k}{B}_kQ_k^T$ and $P_{k+1}{C}_kQ_{k+1}^T$
are the near best rank-$k$ approximation to $A$,
and then become ill-conditioned, because the matrices $B_k$ and $C_k$ have singular values which are rough approximations to the small ones \cite{jia20203, oleary1981}.
Next, therefore, we consider to employ the general-form regularization term
to the projections in (\ref{pro1}) and (\ref{pro2}), which means
we solve the following problem
\begin{equation}\label{mmtik1}
\min_{x\in\mathbb{S}}\|Lx\| \ {\rm subject \ to} \
\mathbb{S}=\left\{x|\left\|P_{k}{B}_kQ_k^Tx-b\right\| = \min\right\}
\end{equation}
or
\begin{equation}\label{mmtik2}
\min_{x\in\mathbb{S}}\|Lx\| \ {\rm subject \ to} \
\mathbb{S}=\left\{x|\left\|P_{k+1}{C}_kQ_{k+1}^Tx-b\right\| = \min\right\},
\end{equation}
with $\ k=1,2,\cdots$.
Problems (\ref{mmtik1}) and (\ref{mmtik2}) can be seen as
projecting the constraint of Problem (\ref{gen1}) onto
a consequence of nested Krylov subspace
and retaining the regularization matrix unchanged.
Therefore, our new algorithms are called the hybrid CGME algorithm and 
hybrid TCGME algorithm,
which are abbreviated as hyb-CGME and hyb-TCGME, respectively.

If $P_{k}{B}_kQ_k^T$ and $P_{k+1}{C}_kQ_{k+1}^T$
are the best rank-$k$ approximation to $A$,
our new algorithms would be the same as the MTSVD method in \cite{hansen92}.
In this sense, the MTSVD method is also a
hybrid method, so is the MTRSVD method \cite{jia2018}.

With regard to the solution $x_{L,k}^{cgme}$ to (\ref{mmtik1})
and $x_{L,k}^{tcgme}$ to (\ref{mmtik2}),
we establish the following results.
\begin{theorem}\label{thm}
Let $Q_{k}$ and $Q_{k+1}$ be defined in (\ref{lsqr1}).
Then the solution to (\ref{mmtik1}) is
\begin{equation}\label{eq5}
x_{L,k}^{cgme} =x_k^{cgme} -\left(L(I_n-Q_kQ_k^T)\right)^{\dagger} Lx_k^{cgme},
\end{equation}
and the solution to (\ref{mmtik2}) can be written as
\begin{equation}\label{eq52}
x_{L,k}^{tcgme} =x_{k}^{tcgme} -
\left(L(I_n-Q_{k+1}Q_{k+1}^T)\right)^{\dagger} Lx_{k}^{tcgme},
\end{equation}
where $x_k^{cgme}$ and $x_{k}^{tcgme}$
are defined in (\ref{xk})
and (\ref{txk}), respectively.
\end{theorem}

{\em Proof}.
Let $P_{cgme}=P_{k}{B}_kQ_{k}^T$ and $P_{tcgme}=P_{k+1}{C}_kQ_{k+1}^T$. From Eld$\acute{e}$n \cite{elden82}, we derive
at iteration $k$, the solution to (\ref{mmtik1}) is
\begin{equation*}
x_{L,k}^{cgme} = (I_n - (L(I_n-P_{cgme}^{\dagger}
P_{cgme}))^{\dagger}L)P_{cgme}^{\dagger}b
\end{equation*}
and to (\ref{mmtik2}) is
\begin{equation*}
x_{L,k}^{tcgme} = (I_n - (L(I_n-P_{tcgme}^{\dagger}
P_{tcgme}))^{\dagger}L)P_{tcgme}^{\dagger}b.
\end{equation*}
Combining the above equations with (\ref{xk})
as well as (\ref{txk}),
we obtain
\begin{equation}\label{eq6}
x_{L,k}^{cgme} = (I_n - (L(I_n-P_{cgme}^{\dagger}
P_{cgme}))^{\dagger}L)x_k^{cgme}.
\end{equation}
and
\begin{equation}\label{eq62}
x_{L,k}^{tcgme} = (I_n - (L(I_n-P_{tcgme}^{\dagger}
P_{tcgme}))^{\dagger}L)x_{k}^{tcgme}.
\end{equation}
Since ${B}_k$ is invertible and $P_{k}$ and $Q_k$
are matrices with orthonormal columns,
it is easy to see that
\begin{equation}\label{eq7}
P_{cgme}^{\dagger}P_{cgme}= Q_kQ_k^T.
\end{equation}
As $C_k$ is the best rank-$k$ approximation to $B_k$, and
$P_{k+1}$ and $Q_{k+1}$ are matrices with orthonormal columns,
we have
\begin{equation}\label{eq8}
P_{tcgme}^{\dagger}P_{tcgme}= Q_{k+1}Q_{k+1}^T.
\end{equation}
Throwing (\ref{eq7}) into (\ref{eq6}) yields (\ref{eq5}) and,
similarly, bringing (\ref{eq8}) into (\ref{eq62}) derives (\ref{eq52}).
\qquad

\begin{remark}
From the form of the regularized solutions (\ref{eq5}) and (\ref{eq52}),
these two solutions can be thought of as iterative solutions (\ref{xk}) and
(\ref{txk}) derived by CGME and TCGME are modified by
$\left(L(I_n-Q_kQ_k^T)\right)^{\dagger} Lx_k^{cgme}$
and $\left(L(I_n-Q_{k+1}Q_{k+1}^T)\right)^{\dagger} Lx_{k}^{tcgme}$, respectively.
\end{remark}
\begin{remark}
When $L=I_n$, recalling (\ref{xk}) and $Q_k^TQ_k=I_n$, it follows from (\ref{eq5}) that
\begin{eqnarray*}
x_{I_n,k}^{cgme} &=&x_k^{cgme} -\left((I_n-Q_kQ_k^T)\right)^{\dagger} x_k^{cgme}\\
&=&(I_n-(I_n-Q_kQ_k^T))x_k^{cgme}\\
&=&Q_kQ_k^Tx_k^{cgme}=Q_kQ_k^TQ_k{B}_k^{-1}P_k^Tb\\
&=&Q_k{B}_k^{-1}P_k^Tb=x_k^{cgme}
\end{eqnarray*}
with $\left((I_n-Q_kQ_k^T)\right)^{\dagger}=I_n-Q_kQ_k^T$
because $(I_n-Q_kQ_k^T)$ is an orthogonal projection.
In the same way, from (\ref{txk}), (\ref{eq52}) and $Q_{k+1}^TQ_{k+1}=I_n$, we derive
\begin{eqnarray*}
x_{I_n,k}^{tcgme} &=&x_{k+1}^{tcgme} -
\left((I_n-Q_{k+1}Q_{k+1}^T)\right)^{\dagger}x_{k}^{tcgme}\\
&=&(I_n-(I_n-Q_{k+1}Q_{k+1}^T)x_{k}^{tcgme}\\
&=&x_{k}^{tcgme}.
\end{eqnarray*}
This means that when $L=I_n$, the regularized solutions obtained by our hybrid algorithms are
the ones to standard-form regularization.
\end{remark}

Next, we consider how to calculate the solutions $x_{L,k}^{cgme}$
and $x_{L,k}^{tcgme}$.
Let $z_k^{cgme} = (L(I_n-Q_kQ_k^T))^{\dagger}Lx_k^{cgme}$
and $z_{k}^{tcgme} = (L(I_n-Q_{k+1}Q_{k+1}^T))^{\dagger}Lx_{k}^{tcgme}$.
Then $z_k^{cgme}$ and $z_{k}^{tcgme}$ are the minimum 2-norm solutions to the least
squares problems
\begin{equation}\label{eq10}
\min_{z\in \mathbb{R}^n}\left\|\left(L(I_n-Q_kQ_k^T)\right)z - Lx_k^{cgme}\right\|
\end{equation}
and
\begin{equation}\label{eq102}
\min_{z\in \mathbb{R}^n}\left\|\left(L(I_n-Q_{k+1}Q_{k+1}^T)\right)z
- Lx_{k}^{tcgme}\right\|,
\end{equation}
respectively.


Evidently, problems (\ref{eq10}) and (\ref{eq102}) have the same structures.
The only difference between them is
$Q_{k+1}$ in (\ref{eq102}), which has one more column than
$Q_k$ in (\ref{eq10}).
As a consequence, in what follows, we
take (\ref{eq10}) as an example to study how to solve them.

Due to the large size of $L(I_n-Q_kQ_k^T)$,
we suppose that the problems (\ref{eq10})
can only be solved by iterative algorithms.
We will use the LSQR algorithm \cite{lsqr1982}
to solve the problem. To take fully advantage of
the sparsity of $L$ itself and reduce the overhead of computation
and the storage memory, it is critical to avoid forming the dense matrix
$L(I_n-Q_kQ_k^T)$ explicitly within LSQR.
Note that the only action of $L(I_n-Q_kQ_k^T)$ in the Lanczos diagonalization process is to form the products
of it and its transpose with vectors. We propose Algorithm~\ref{alg2},
which efficiently implements the Lanczos bidiagonalization process
without forming $L(I_n-Q_kQ_k^T)$ explicitly.

\begin{algorithm}
\caption{$\widehat k$-step Lanczos bidiagonalization process on $L(I_n -Q_kQ_k^T)$.}
\label{alg2}
\begin{algorithmic}[1]
\State Taking $\beta_1\widehat{u}_1= Lx_k$, $\beta_1=\|Lx_k\|$, $w_1= L^T\widehat{u}_1$, and define $\widehat{\beta}_1\widehat{v}_0 = 0$.
\For{$j=1, 2, \ldots, \widehat k$} 
\State $\widehat{r}=w_j-Q_k(Q_k^Tw_j)-\widehat{\beta}_j\widehat{v}_{j-1}$
\State $\widehat{\alpha}_j=\|\widehat{r}\|$;~~$\widehat{v}_j=\widehat{r}/\widehat{\alpha}_j$; ~$g_{j} = Q_k^T\widehat{v}_{j}$
\State $\widehat{s}=L\widehat{v}_j-L(Q_kg_j)-\widehat{\alpha}_j\widehat{u}_j$
\State $\widehat{\beta}_{j+1} = \|\widehat{s}\|$;~~$\widehat{u}_{j+1} = \widehat{s}/\widehat{\beta}_{j+1}$;~$w_{j+1} = L^T\widehat{u}_{j+1}$
\EndFor
\end{algorithmic}
\end{algorithm}

We now consider the solution of (\ref{eq10}) using LSQR.
Let
\begin{equation}\label{Q}
Q=\left(
\begin{array}{cc}
Q_k & Q_k^{\perp} \\
\end{array}
\right)
\end{equation}
be an orthogonal matrix and $Q_k^{\perp}\in\mathbb{R}^{n\times (n-k)}$
an orthogonal complement of the matrix $Q_k$.
With the notation of (\ref{Q}), we have
\begin{equation}\label{LQ}
L(I_n -Q_kQ_k^T)=LQ_k^{\perp}(Q_k^{\perp})^T.
\end{equation}
Because $Q_k^{\perp} $ is column orthonormal,
the nonzero singular values of $LQ_k^{\perp}(Q_k^{\perp})^T$
are identical to the singular values of $LQ_k^{\perp}$,
we have
\begin{equation}\label{condi}
\kappa(L(I_n-Q_kQ_k^T))=\kappa(LQ_k^{\perp}(Q_k^{\perp})^T)
=\kappa(LQ_k^{\perp}).
\end{equation}

We now investigate how the conditioning of (\ref{eq10}) changes as $k$ increases with notation.
\begin{theorem}\label{thm3}
Let the matrix $Q_k^{\perp}$ be defined in (\ref{Q}) and $L\in\mathbb{R}^{p\times n}$.
When $p\geq n-k$, we obtain
\begin{equation}\label{eq3}
\kappa(LQ_k^{\perp})\geq\kappa(LQ_{k+1}^{\perp}),\quad k = 2, 3, \cdots, n-1,
\end{equation}
that is,
\begin{equation}\label{eq2}
\kappa(L(I_n -Q_kQ_k^T))\geq\kappa(L(I_n -Q_{k+1}Q_{k+1}^T)),
\quad k = 2, 3, \cdots,n-1.
\end{equation}
\end{theorem}

{\em Proof}.
Exploiting the lemma presenting in \cite[pp.78]{golub13} yields
\begin{eqnarray*}
\sigma_{\max}(LQ_k^{\perp})\geq \sigma_{\max}(LQ_{k+1}^{\perp}),\label{1}\\
\sigma_{\min}(LQ_k^{\perp})\leq \sigma_{\min}(LQ_{k+1}^{\perp}).\label{2}
\end{eqnarray*}
From the above of the inequalities,
we derive
\begin{equation*}
\kappa(LQ_k^{\perp})=\frac{\sigma_{\max}(LQ_k^{\perp})}{\sigma_{\min}(LQ_k^{\perp})}
\geq \frac{\sigma_{\max}(LQ_{k+1}^{\perp})}{\sigma_{\min}(LQ_{k+1}^{\perp})}
=\kappa(LQ_{k+1}^{\perp}),
\end{equation*}
that is, (\ref{eq3}) holds. Noticing (\ref{LQ}),
we directly obtain (\ref{eq2}) from (\ref{eq3}).
\qquad

Theorem~\ref{thm3} indicates that, when applied to solve (\ref{eq10}),
the LSQR algorithm generally converges faster with $k$
because the worst LSQR convergence factor is
$\frac{\kappa(LQ_k^{\perp})-1}{\kappa(LQ_k^{\perp})+1}$,
see \cite[p.~291]{bjorck96}.
In particular, in exact arithmetic,
LSQR will find the exact solution of (\ref{eq10})
after at most $n-k$ iterations.
Exploiting the same method and technology,
we can derive similar results for (\ref{eq102}).

Now we present the hyb-CGME and hyb-TCGME algorithms,
named Algorithm~\ref{alg:pcgme} and Algorithm~\ref{alg:ptcgme}, respectively.

\begin{algorithm}
\caption{(hyb-CGME) \ Given $A\in \mathbb{R}^{m\times n}$,
$L\in\mathbb{R}^{p\times n}$ and $b\in\mathbb{R}^n$, compute the solution
$x_{L,k}^{cgme}$.}
\label{alg:pcgme}
\begin{algorithmic}[1]
\State Use Algorithm \ref{alg1} to compute the projection $P_{k}^TB_kQ_k$.
\State Compute $x_k^{cgme}$ by (\ref{xk}).
\State Compute the solution $z_k^{cgme}$ to (\ref{eq10}) by LSQR with Algorithm~\ref{alg2}.
\State Compute the solution $x_{L,k}^{cgme}=x_k^{cgme}-z_k^{cgme}$ .
\end{algorithmic}
\end{algorithm}

\begin{algorithm}
\caption{(hyb-TCGME) \ Given $A\in \mathbb{R}^{m\times n}$,
$L\in\mathbb{R}^{p\times n}$ and $b\in\mathbb{R}^n$, compute the solution
$x_{L,k}^{tcgme}$.}\label{alg:ptcgme}
\begin{algorithmic}[1]
\State Use Algorithm~\ref{alg1} to compute the projection $P_{k+1}^T{C}_kQ_{k+1}$.
\State Compute $x_{k}^{tcgme}$ by (\ref{txk}).
\State Compute the solution $z_{k}^{tcgme}$ to (\ref{eq102}) by LSQR with Algorithm~\ref{alg2} replacing $Q_k$ with $Q_{k+1}$.
\State Compute the solution $x_{L,k}^{tcgme}=x_{k}^{tcgme}-z_{k}^{tcgme}$ .
\end{algorithmic}
\end{algorithm}

We emphasize that, in step 3 of Algorithm~\ref{alg:pcgme}-\ref{alg:ptcgme},
we will use LSQR with Algorithm~\ref{alg2} to solve (\ref{eq10}) and (\ref{eq102}) with a given tolerance $tol$ as the
stopping criterion that substantially exploits the Matlab built-in function \textsf{lsqr.m}
without explicitly computing
the matrix products $L(I_n-Q_kQ_k^T)$ and $L(I_n-Q_{k+1}Q_{k+1}^T)$.
We will make an insightful analysis in next section and
show that the default $tol=10^{-6}$ is generally
good enough and larger $tol$ can be allowed in practical applications.

\section{The stopping tolerance $tol$ of LSQR for solving the inner least squares}\label{sec4}

In this section we establish some important theoretical results on our proposed new algorithms and get insight into how to determine the stopping tolerance of LSQR which is used to solve the
inner least squares problems.
We only take the hyb-CGME algorithm as an example to analyze the accuracy of the calculated and accurate regularized solution because the hyb-TCGME algorithm is the similar with it.
For convenience of writing and without any ambiguity,
in what follows, we drop the superscript,
which means we will use $z_k$,
$x_{k}$ and $x_{L,k}$ instead of $z_k^{cgme}$, $x_{k}^{cgme}$ and $x_{L,k}^{cgme}$, respectively.

First of all, with the notation in Section \ref{sec3},
we establish the result of the accuracy of the computed solution $\bar{z}_k$
with the stopping tolerance $tol$ in the following lemma.
\begin{lemma}
Let $z_k$ and $\bar{z}_k$ be the exact solution and
the computed solution by LSQR with the stopping criterion $tol$
for the problem (\ref{eq10}), respectively.
Then
\begin{equation}\label{relerr}
\frac{\|z_k-\bar{z}_k\|}{\|z_k\|}\leq
\frac{\kappa(LQ_k^{\perp})}
{1-tol\cdot
\kappa(LQ_k^{\perp})}\left(2+\frac{\kappa(LQ_k^{\perp})\|r\|}
{\|LQ_k^{\perp}\|\|z_k\|}
\right)tol.
\end{equation}
with $r= Lx_k-L(I_n -Q_kQ_k^T)z_k$.
\end{lemma}

{\em Proof}.
From \cite{lsqr1982} it follows that
with the stopping tolerance $tol$
the computed solution $\bar{z}_k$ is the exact
solution to the perturbed problem
\begin{equation}\label{perb}
\min_{z\in \mathbb{R}^{n}}\|(L(I_n -Q_kQ_k^T)+E_k)z - Lx_k\|,
\end{equation}
where $$
E_k=-\frac{r_kr_k^TL(I_n -Q_kQ_k^T)}{\|r_k\|^2}
$$
is the perturbation matrix
with $
r_k=Lx_k
-L(I_n -Q_kQ_k^T)\bar{z}_k
$
and
$$
\frac{\|E_k\|}{\|L(I_n -Q_kQ_k^T)\|}
=\frac{\|(I_n -Q_kQ_k^T)L^Tr_k\|}
{\|L(I_n -Q_kQ_k^T)\|\|r_k\|} \leq tol.
$$
Exploiting
the standard perturbation theory \cite[p.~382]{higham02} and (\ref{condi}), we obtain
\begin{equation}\label{relerror}
\frac{\|z_k-\bar{z}_k\|}{\|z_k\|}\leq
\frac{\kappa(LQ_k^{\perp})}
{1-tol\cdot
\kappa(LQ_k^{\perp})}\left(2+(\kappa(LQ_k^{\perp})+1)\frac{\|r\|}
{\|LQ_k^{\perp}\|\|z_k\|}
\right)tol,
\end{equation}
According to the proof of Wedin's Theorem; see \cite[pp.~400]{higham02},
we find that the factor
$\kappa(LQ_k^{\perp})+1$ in (\ref{relerror})
can be replaced by $\kappa(LQ_k^{\perp})$
in the context
since the right-hand side $Lx_k$ in (\ref{perb})
is unperturbed.
Therefore, it is easy to get (\ref{relerr})
from (\ref{relerror}).
\qquad

In applications, $L$ is typically well conditioned
\cite{hansen98,hansen10}, so is $LQ_k^{\perp}$ because
of the orthonormality of $Q_k^{\perp}$. Therefore,
the left hand side of (\ref{relerr}) is at least as small as
${\cal O}(tol)$ with a generic constant in ${\cal O}(\cdot)$,
which means
\begin{equation}\label{appz_k}
\frac{\|z_k-\bar{z}_k\|}{\|z_k\|} \leq {\cal O}(tol).
\end{equation}


Keep in mind that the hyb-CGME solution can be written as
$x_{L,k}=x_k-z_k.$ Define the computed hyb-CGME solution
$$\bar{x}_{L,k}=x_k-\bar{z}_k.$$
Therefore, we derive
\begin{equation}\label{lkk}
\|x_{L,k}-\bar{x}_{L,k}\|=\|z_k-\bar{z}_k\|.
\end{equation}
From (\ref{appz_k}) and (\ref{lkk})
it is reasonable to suppose
\begin{equation}\label{xkzk}
\frac{\|x_{L,k}-\bar{x}_{L,k}\|}{\|x_{L,k}\|}
\leq {\cal O}(tol)
\end{equation}
since it is generally impossible that
$\|x_{L,k}\|$ is much smaller or larger than $\|z_k\|$.

Let $x_{L}^{opt}$ be an {optimal} regularized solution to
the problem (\ref{gen1}) with the white noise $e$.
Then under a certain necessary discrete
Picard condition, a GSVD analysis indicates
that the error $\|x_{L}^{opt}-x_{true}\|\geq {\cal O}(\|e\|)$ with a
generic constant in ${\cal O}(\cdot)$; see \cite[p.~83]{hansen98} and \cite{jia2018}. 

Define $x_{L,k_0}$ and $\bar{x}_{L,k_0}$ to be the best regularized solution
and the computed best regularized solution by hyb-CGME algorithm, respectively.
Obviously,
\begin{equation}\label{compare}
\|x_{L,k_0}-x_{true}\|\geq\|x_L^{opt}-x_{true}\|\geq {\cal O}(\|e\|).
\end{equation}
With these notations, the following theorem compares the accuracy of
the best regularized solution with the computed best regularized solution.
\begin{theorem}\label{thm6}
	If $L$ is well conditioned and $\|A\|\approx 1$,
	then
	\begin{equation}\label{newadd}
	\frac{\|x_{L,k_0}-x_{true}\|}
	{\|x_{true}\|}\geq {\cal O}\left(\frac{\|e\|}{\|b_{true}\|}\right).
	\end{equation}
	If
	\begin{equation}\label{accucond}
	{\cal O}(tol)< \frac{\|e\|}{\|b_{true}\|},
	\end{equation}
	then
	\begin{equation}\label{computerror}
	\left|\frac{\|\bar{x}_{L,k_0}-x_{true}\|}{\|x_{true}\|}-\frac{\|x_{L,k_0}-x_{true}\|}
	{\|x_{true}\|}\right|\leq {\cal O}(tol)
	\end{equation}
	with a generic constant in ${\cal O}(\cdot)$.
\end{theorem}

{\em Proof}.
Supposing that $Ax_{true}=b_{true}$ which means that the noise-free problem of (\ref{eq1}) is consistent,
we have $\|b_{true}\| \leq \|A\|\|x_{true}\|$.
Using the inequality, it follows from (\ref{compare}) that
$$
\frac{\|x_{L,k_0}-x_{true}\|}
{\|x_{true}\|}\approx \|A\|\frac{\|x_{L,k_0}-x_{true}\|}
{\|b_{true}\|}\geq \|A\|{\cal O}\left(\frac{\|e\|}{\|b_{true}\|}\right).
$$
Because by suitable scaling $\|A\|\approx 1$ can always be done,
we derive (\ref{newadd}) from the above relations.
Using (\ref{relerror}) and (\ref{xkzk}) as well as $\|x_{L,k_0}\|\approx \|x_{true}\|$ obtains
\begin{eqnarray}
\frac{\|\bar{x}_{L,k_0}-x_{true}\|}{\|x_{true}\|}&\leq &\frac{\|x_{L,k_0}-x_{true}\|}
{\|x_{true}\|}+\frac{\|x_{L,k_0}-\bar{x}_{L,k_0}\|}{\|x_{true}\|} \nonumber\\
&=&\frac{\|x_{L,k_0}-x_{true}\|}
{\|x_{true}\|}+\frac{\|x_{L,k_0}-\bar{x}_{L,k_0}\|}{\|x_{L,k_0}\|}
\frac{\|x_{L,k_0}\|}{\|x_{true}\|} \nonumber\\
&\leq&\frac{\|x_{L,k_0}-x_{true}\|}
{\|x_{true}\|}+{\cal O}(tol). \label{oneside}
\end{eqnarray}
On the other hand, we similarly obtain
\begin{equation}\label{twoside}
\frac{\|\bar{x}_{L,k_0}-x_{true}\|}{\|x_{true}\|}\geq \frac{\|x_{L,k_0}-x_{true}\|}
{\|x_{true}\|}-{\cal O}(tol).
\end{equation}
Combining (\ref{oneside}) with (\ref{twoside}) obtains
\begin{equation*}\label{appeq}
	\frac{\|x_{L,k_0}-x_{true}\|}
	{\|x_{true}\|}-{\cal O}(tol)
	\leq \frac{\|\bar{x}_{L,k_0}-x_{true}\|}{\|x_{true}\|}\leq
	\frac{\|x_{L,k_0}-x_{true}\|}{\|x_{true}\|}+{\cal O}(tol),
	\end{equation*}
which implies (\ref{computerror}).
\qquad

Clearly, (\ref{newadd}) and (\ref{computerror})
imply that, under condition (\ref{accucond}),
the computed $\bar{x}_{L,k_0}$ is almost the same
as the exact $x_{L,k_0}$ as an approximation to $x_{true}$.
Furthermore, according to the analysis above,
we can establish more general results, including Theorem~\ref{thm6} as a special case.

\begin{theorem}\label{thm7}
If $L$ is well conditioned and $\|A\|\approx 1$,
then
\begin{equation}\label{newadd2}
\frac{\|x_{L,k}-x_{true}\|}
{\|x_{true}\|}\geq {\cal O}\left(\frac{\|e\|}{\|b_{true}\|}\right),\
k=1,2,\ldots,\min\{m,n\}.
\end{equation}
Let $\bar{x}_{L,k}$ be the computed solutions obtained by
Algorithm~\ref{alg:pcgme} using LSQR in step 3 with the stopping criterion $tol$. Then if
the condition (\ref{accucond}) is satisfied, for $k=1,2,\ldots,k_0$ and a few
$k>k_0$ we have
\begin{equation}\label{computerror2}
	\left|\frac{\|\bar{x}_{L,k}-x_{true}\|}{\|x_{true}\|}-\frac{\|x_{L,k}-x_{true}\|}
	{\|x_{true}\|}\right|\leq {\cal O}(tol)
\end{equation}
with a generic constant in ${\cal O}(\cdot)$.
\end{theorem}

{\em Proof}.
Notice that $x_{L,k_0}$ is the best possible regularized solutions
by hyb-CGME and Algorithm \ref{alg1} is run for $k=1,2,\ldots,\min\{m,n\}$,
i.e.,
\begin{equation}\label{best2}
\frac{\|x_{L,k_0}-x_{true}\|}
{\|x_{true}\|}=\min_{k=1,2,\ldots,\min\{m,n\}} \frac{\|x_{L,k}-x_{true}\|}
{\|x_{true}\|}.
\end{equation}
Because the iterates $\|x_{L,k}\|$
exhibits increasing tendency, and it first approximates $x_{true}$ from below for
$k = 1, 2, \cdots, k_0$,
in later stages, however,
the iterates will start to diverge from
$x_{true}$.
But it won't deviate from $x_{true}$ too much for a few $k > k_0$,
therefore, for $k = 1, 2,..., k_0 \ {\rm and\ a \ few}\ k > k_0$, we have
\begin{equation}\label{estimate}
\frac{\|x_{L,k}\|}{\|x_{true}\|}={\cal O}(1).
\end{equation}
From (\ref{best2}), (\ref{estimate})
and the proof of Theorem~\ref{thm6},
we attain that when the index $k_0$ is replaced
by $k=1,2,\ldots,k_0$ and {a few} $k>k_0$,
(\ref{newadd}) and
(\ref{computerror}) also hold under the condition (\ref{accucond}).
\qquad

Clearly, (\ref{newadd2}) and (\ref{computerror2})
imply that, under condition (\ref{accucond}),
the calculated $\bar{x}_{L,k}$ is almost the same
as the exact $x_{L,k}$ as an approximation to $x_{true}$
for $k=1,2,\ldots,k_0$ and a few $k>k_0$.
Although the above analysis is for using LSQR to solve (\ref{eq10}),
taking advantage of the same method and technique,
for solving (\ref{eq102}) we can get exactly the same results.

What we have to highlight is that
the relative noise level $\frac{\|e\|}{\|b_{true}\|}$ is typically
more or less around $10^{-3}$ or $10^{-2}$ in practice applications,
three or four orders bigger than
$10^{-6}$.
From the above analysis,
we draw the conclusion that it is generally enough to set
$tol=10^{-6}$ in LSQR at step 3 of Algorithms~\ref{alg:pcgme}--\ref{alg:ptcgme}.
A smaller $tol$ will result in more inner
iterations without any gain in the accuracy of $\bar{x}_{L,k}$
for $k=1,2,\ldots,k_0$ and a few $k>k_0$.
Moreover, Theorems~\ref{thm6}--\ref{thm7} indicate that $tol=10^{-6}$
is generally well conservative and larger $tol$ can be used, so
that LSQR fewer iterations to achieve the convergence
and the hyb-CGME algorithm and the hyb-TCGME algorithm
are more efficient.

In summary, the conclusion is that a widely varying choice of $tol$ has no
effects on the regularized solution of hyb-CGME and hyb-TCGME, provided
that $tol<\frac{\|e\|}{\|b_{true}\|}$ is considerably and
the regularization matrix $L$ is well conditioned, but has substantial effects on the efficiency of hyb-CGME and hyb-TCGME.

Finally, there is something we should say on the efficiency comparison of Algorithm~\ref{alg:pcgme}-\ref{alg:ptcgme} with the JBDQR algorithm
which is more efficient than the algorithm in \cite{kilmer07}; see \cite{jia2020} for more details.
It is obvious that at iteration $k$,
Algorithm \ref{alg:pcgme} and \ref{alg:ptcgme} cost almost the same computation, since
the truncation in Algorithm \ref{alg:ptcgme} costs little.
And they may be considerably cheaper than JBDQR because JBDQR has to compute a large-scale least squares problem with coefficient matrix
$(A^T,L^T)^T$ which is larger than the problem itself and may be very expensive.
Compared with JBDQR, Algorithm \ref{alg:pcgme} and Algorithm \ref{alg:ptcgme} only need to solve a large-scale least squares problem with coefficient matrices $L(I_n -Q_kQ_k^T)$ and $L(I_n -Q_{k+1}Q_{k+1}^T)$, respectively,
and it turns out that they become better conditioned with $k$.
This means used to solve the
inner least squares problems in hyb-CGEM and hyb-TCGME, LSQR converges faster as $k$ increases.
The faster LSQR converges,
the cheaper our proposed new hybrid algorithms.
In addition, to obtain the regularized solution, JBDQR finally has to solve a large-scale linear
least squares problem with coefficient matrix
$(A^T,L^T)^T$, which may also be particularly expensive.

\section{Numerical experiments}
In this section, we report numerical experiments
on several problems, including one dimensional and two dimensional problems,
to demonstrate that the proposed new hybrid algorithms
work well and that the regularized solutions obtained by hyb-TCGME are at least as accurate
as those obtained by JBDQR,
and the proposed new method is considerably more efficient than JBDQR.
All the computations are carried out in Matlab R2019a 64-bit on
11th Gen Intel(R) Core(TM) i5-1135G7 2.40GHz processor
and 16.0 GB RAM with the machine precision
$\epsilon_{\rm mach}= 2.22\times10^{-16}$ under the Miscrosoft
Windows 10 64-bit system.

Table~\ref{tab1} lists all test problems, each problem using its default parameter(s).
These problems are one- and two-dimensional problems from the
regularization toolbox \cite{hansen2007(2),berisha2012restore,nagy2004} and
the Matlab Image Processing Toolbox. 
The one-dimensional problems shaw and baart are severely ill-posed, 
while heat is moderately ill-posed, 
and deriv2 is mildly ill-posed, all of which are of the order $m=n=10,000$.
The two-dimensional problem mri is from the Matlab Image Processing Toolbox
with an order of $m=n=16,284$. 
Other two-dimensional problems are from \cite{nagy2004} with an order 
of $m=n=65, 536$. We denote the relative noise level
$$
\varepsilon = \frac{\|e\|}{\|b_{true}\|}.
$$
\begin{table}[ht]
\caption{The description of test problems.}\label{tab1}
\begin{tabular}{@{}lll@{}}
\toprule
 Problem        & Description                                & Size of $m, n$ \\
\midrule
      {shaw}     & one dimensional image restoration model    & $m=n=10,000$\\
     {baart}    & one dimensional gravity surveying problem  &  $m=n=10,000$\\
     {heat}     & Inverse heat equation                   &  $m=n=10,000$\\
     {deriv2}   & Computation of second derivative        &  $m=n=10,000$\\
     {mri}      & Two dimensional image deblurring      & $m=n=16,284$\\
     {grain}    & spatially variant Gaussian blur      & $m=n=65, 536$\\
     {satellite} &spatially invariant atmospheric turbulence        & $m=n=65, 536$ \\
     {GaussianBlur440} & spatially invariant Gaussian blu      & $m=n=65, 536$\\
\bottomrule
\end{tabular}
\end{table}

Let $x^{reg}$ denote the regularized solution obtained by
algorithms. We use the relative error
\begin{equation}\label{rel}
\frac{\|L(x^{reg}-x_{true})\|}{\|Lx_{true}\|}
\end{equation}
to plot the convergence curve of each algorithm with respect to $k$,
which is more instructive and suitable to use
the relative error (\ref{rel}) in the general-form regularization
context other than the standard relative error of $x_{reg}$;
see \cite[Theorems 4.5.1-2]{hansen98} and \cite{jia2020} for more details.
In the tables to be presented,
we will list the smallest relative errors and the total outer iterations
required to obtain the smallest relative errors in the braces.
We also will list the total CPU time which is counted in seconds
by the Matlab built-in commands
{\sf tic} and {\sf toc}
and the corresponding total outer iterations.
For the sake of length, 
we only display the noise levels $\varepsilon =10^{-1}, \ 5 \times10^{-2},\ {\rm and} \ 10^{-2}$ in Tables \ref{tab2} and \ref{tab3}.

For our new algorithms hyb-CGME and hyb-TCGME,
we use the Matlab built-in function {\sf lsqr}
with Algorithm~\ref{alg2} avoiding forming the dense matrices $L(I_n-Q_kQ_k^T)$
and $L(I_n-Q_{k+1}Q_{k+1}^T)$ explicitly
to compute (\ref{eq10}) and (\ref{eq102})
with the default stopping tolerance $tol=10^{-6}$.
For the JBDQR algorithm \cite{jia2020},
we use the same function with the same stopping tolerance
to solve
the inner least squares problems.

For the four one-dimensional test problems we use the code of \cite{hansen2007(2)} to generate $A$,
the true solution $x_{true}$ and noise-free right-hand side $b_{true}$.
Purely for test purposes, we choose $L=L_1$ defined by
\begin{equation}\label{l1}
L_1 = \left(
        \begin{array}{ccccc}
          1 & -1 &  &  &  \\
           & 1 & -1 &  &  \\
           &  & \ddots & \ddots &  \\
             &  &  & 1  & -1\\
        \end{array}
      \right)\in \mathbb{R}^{(n-1)\times n}.
\end{equation}
For the two-dimensional image deblurring problems,
the regularization matrix is chosen as
\begin{equation}\label{l3}
L=\left(
    \begin{array}{c}
      I_N\otimes L_1 \\
      L_1 \otimes I_N \\
    \end{array}
  \right),
\end{equation}
where $L_1$ is defined in (\ref{l1}),
which is the scaled discrete approximation of the first
derivative operator in the two dimensional case incorporating no assumptions
on boundary conditions; see \cite[Chapter 8.1-2]{hansen10}.

\begin{table}[ht]
\caption{The relative errors and the optimal
regularization parameters in the braces for test problems in Table~\ref{tab1}.}\label{tab2}
\begin{minipage}[t]{1\textwidth}
\begin{tabular*}{\textwidth}{@{\extracolsep\fill}lcccccc}
\toprule%
& \multicolumn{3}{@{}c@{}}{$\varepsilon=10^{-1}$}  \\\cmidrule{2-4}%
&{hyb-CGME}    &JBDQR    &{hyb-TCGME} \\
\midrule
{shaw}     &0.9908(2)     &3.2943(1)      &0.2244(7)  \\
{baart}    &0.9663(2)      &0.6093(1)       &0.5615(4)    \\
{heat}     &0.9516(4)      &0.3636(1)       &0.3689(13)  \\
{deriv2}   &0.8002(2)      &1.2761(1)       &0.4805(6)   \\
{mri}       &0.9883(2)     &0.9134(6)      &0.9169(12)   \\
{grain}     &0.9992(4)     &0.9524(9)      &0.9601(31)   \\
{sate} &0.9968(6)     &0.9681(13)     &0.9682(47)   \\
{blur440}   &0.9913(2)     &0.9711(6)      &0.9721(21)  \\
\bottomrule
\end{tabular*}
\end{minipage}
\begin{minipage}[t]{1\textwidth}
\begin{tabular*}{\textwidth}{@{\extracolsep\fill}lcccccc}
& \multicolumn{3}{@{}c@{}}{$\varepsilon=5\times10^{ -2}$} \\\cmidrule{2-4}%
&{hyb-CGME}    &JBDQR    &{hyb-TCGME}  \\
\midrule
{shaw}    &0.9770(3)     &1.3390(1)     &0.2515(7)     \\
{baart}   &0.9869(2) 	 &0.6063(1)	    &0.5535(3)	   \\
{heat}    &0.7805(5)	 &0.3470(1)	    &0.3499(15)	    \\
{deriv2}  &0.5907(2)	 &0.9265(1)	    &0.4443(9)	   \\
{mri}        &0.9540(4)	    &0.8840(12)	    &0.8873(20)	    \\
{grain}      &0.9975(7) 	&0.9070(19)	    &0.9237(43)	   \\
{sate}  &0.9924(11)	&0.9597(26)	    &0.9597(61)	   \\
{blur440}    &0.9855(4)     &0.9654(13)     &0.9666(39)     \\
\bottomrule
\end{tabular*}
\end{minipage}
\begin{minipage}[t]{1\textwidth}
\begin{tabular*}{\textwidth}{@{\extracolsep\fill}lcccccc}
& \multicolumn{3}{@{}c@{}}{$\varepsilon=10^{-2}$ } \\\cmidrule{2-4}%
&{hyb-CGME}    &JBDQR    &{hyb-TCGME} \\
\midrule
{shaw}     &0.9681(4)      &0.3237(1)       &0.1972(7)     \\
{baart}    &0.8803(3)      &0.7054(1)       &0.5500(3)    \\
{heat}     &0.5695(9)      &0.2146(5)       &0.2128(20)   \\
{deriv2}   &0.6113(2)      &0.4359(1)       &0.6625(3)   \\
{mri}         &0.8838(12)      &0.8412(51)       &0.8448(68)   \\
{grain}       &0.9835(22)      &0.7089(71)       &0.8005(57)    \\
{sate}   &0.9713(35)      &0.9321(132)      &0.9345(117)   \\
{blur440}     &0.9717(12)      &0.9511(69)       &0.9528(159)     \\
\bottomrule
\end{tabular*}
\end{minipage}
\end{table}

From Table \ref{tab2},
we observe that for all test problems,
the best regularized solution by {hyb-TCGME} is
at least as accurate as and can be considerably more accurate than
that by JBDQR;
see, e.g., the results on test problems shaw and baart.
And the regularized solution by {hyb-CGME} is less
accurate than that by JBDQR and {hyb-TCGME},
except for the test problem shaw with $\varepsilon=10^{-1}$
and $\varepsilon=5\times 10^{-2}$.
Furthermore, for the test problems{ shaw} and {deriv2},
we can see from the table that 
the JBDQR algorithm fails
while {hyb-CGME} and {hyb-TCGME} work very well,
especially {hyb-TCGME} which derives the best regularized solutions with high accuracy
with noise levels $\varepsilon=10^{-1}$ and $5\times10^{-2}$.
For two-dimensional test problems, 
we can observe that the best regularized solutions by our new hybrid
has almost the same accuracy with the ones by JBDQR.

Insightfully, we observe from comparing the second column
with the fourth column of Table \ref{tab2} that for each test
problem and given $\varepsilon$,
hyb-TCGME obtains the best regularized solution that is correspondingly more accurate
and requires a correspondingly larger regularization parameter than hyb-CGME.
This indicates that, compared with hyb-TCGME,
hyb-CGME captures less dominant GSVD components of $\{A, L\}$.
Recall the comments at the end of Section \ref{sec2}.
CGME only has partial regularization and the rank-$k$ approximation to $A$
obtained by it is less accurate than the counterpart by TCGME.
Therefore, the less accurate rank-$k$ approximation to $A$ generated by the projection,
the less accurate of the generalized solution obtained by our new hybrid method. 
This is reasonable because a
good regularized solution must capture all the needed dominant GSVD components
of the matrix pair $\{A, L\}$ and, meanwhile, suppress those corresponding to small
generalized singular values;
see, e.g.,\cite{hansen98,hansen10,kilmer07,jia2020}.


\begin{table}[ht]
\caption{The CPU time of {hyb-CGME}, JBDQR and {hyb-TCGME},
the ratio of the CPU time of JBDQR to that of hyb-CGME  ({times1}),
the ratio of the CPU time of JBDQR to that of hyb-TCGME ({times2}), and
the total outer iterations ({iteration})
for test problems in Table~\ref{tab1}.}\label{tab3}
\begin{minipage}[t]{1\textwidth}
\begin{tabular*}{\textwidth}{@{\extracolsep\fill}lcccccc}
\toprule%
&\multicolumn{6}{@{}c@{}}{$\varepsilon=10^{-1}$} \\\cmidrule{2-7}%
&{hyb-CGME}    &JBDQR    &{hyb-TCGME}  &{ times1} &{ times2} &{ iteration}\\
\midrule
{shaw}    &6.5319  &2678.8398 &7.4490 &410.1164	&359.6241   &12\\
{baart}    &3.0381    &1940.6978       &3.0063    &638.7866	&645.5436 &12  \\
{heat}     &11.9812    &2744.7590       &8.5197    &229.0888	&322.1662 &20 \\
{deriv2}   &20.1345     &8718.1975      &26.2901   &432.9979	&331.6152  &20  \\
{grain} &33.9623  &339.2248   &37.3480     &9.9883	&9.0828  &50        \\
{mri}    &7.1192    &47.9356      &6.9455    &6.7333	&6.9017 &50  \\
{satellite}   &25.4080    &230.6816      &36.3500    &9.0791	&6.3461 &50 \\
{blur440}   &24.8721      &57.2304      &42.8874    &2.3010	 &1.3344 &50  \\
\bottomrule
\end{tabular*}
\end{minipage}
\begin{minipage}[t]{1\textwidth}
\begin{tabular*}{\textwidth}{@{\extracolsep\fill}lcccccc}
& \multicolumn{6}{@{}c@{}}{$\varepsilon=5\times10^{ -2}$} \\\cmidrule{2-7}%
&{hyb-CGME}    &JBDQR    &{hyb-TCGME} &{ times1} &{ times2} &{ iteration} \\
\midrule
{shaw}    &7.7211    &2652.1893     &8.2861  &343.4989	&320.0769 &12     \\
{baart}   &3.7360	 &1815.0300	    &3.2283  &485.8217  &562.2247 &12	   \\
{heat}    &13.2105	 &2985.5863     &8.6991  &226.0010  &343.2063 &20	    \\
{deriv2}  &10.6474	 &4719.4951     &13.7445 &443.2533	&343.3733 &18   \\
{mri}    &26.4831    &120.9110      &26.4821    &4.5656	&4.5658 &100  \\
{grain}    &76.2086     &730.3901    &72.8844     &9.5841	&10.0212 &100 \\
{satellite}  &62.3031    &443.2401  &71.8914   &7.1143	&6.1655 &100 \\
{blur440}   &73.8359     &112.0466  &46.0917 &1.5175	&2.4309 &100  \\
\bottomrule
\end{tabular*}
\end{minipage}
\begin{minipage}[t]{1\textwidth}
\begin{tabular*}{\textwidth}{@{\extracolsep\fill}lcccccc}
& \multicolumn{6}{@{}c@{}}{$\varepsilon=10^{-2}$ } \\\cmidrule{2-7}%
&{\ hyb-CGME}    &JBDQR    &{hyb-TCGME} &{ times1} &{ times2} &{ iteration}\\
\midrule
{shaw}     &8.3601    &4120.1932   &8.7147    &492.8401	&472.7865 &16        \\
{baart}    &3.1327    &1966.6141   &3.0570    &627.7697	&643.3150 &12  \\
{heat}     &15.4444   &7149.1920   &12.5738   &462.8987	&568.5785 &30 \\
{deriv2}   &11.6898   &4937.2102   &14.9254   &422.3520 &330.7925 & 20 \\
{grain}     &365.3413	&3094.4985	&249.2137&8.4702	&12.4170
&300   \\
{mri}    &435.0728	&1976.4200	&422.5347 &4.5427	&4.6775
&500     \\
{satellite}  &281.9240      &1764.3952       &301.2180 &6.2584	&5.8575
&400     \\
{blur440}  &917.5464	&1975.2774	&795.1023 &2.1528	&2.4843
& 1000 \\
\bottomrule
\end{tabular*}
\end{minipage}
\end{table}

Recall the analysis of the method based on JBD process of Section \ref{sec2} and the comments of Section \ref{sec4}.
At the same outer iterations, hyb-CGME and hyb-TCGME may be much cheaper than JBDQR,
and they may cost the same overhead because the truncation in hyb-TCGME is very cheap.
As can be seen from Table \ref{tab3},
for each test problem, the CPU time of {hyb-CGME} and
{hyb-TCGME} is
significantly less than that of JBDQR at the same outer iterations.
For the severely ill-posed problems shaw and baart,
hyb-CGME and hyb-TCGME take less than ten seconds, respectively,
while JBDQR takes almost two thousands of seconds;
in particular, for the test problem baart with noise level of $\varepsilon=10^{-2}$,
our new method only takes less than 13 seconds,
whereas JBDQR takes more than even thousands of seconds.
For the moderately and mildly ill-posed problems heat and deriv2,
the CPU time of hyb-CGME and hyb-TCGME is less than 20 seconds except
for deriv2 with the noise level $\varepsilon=10^{-1}$,
the CPU time of hyb-TCGME is 26.2901, a little bit more than 20 seconds.
In addition, we also observe from the table that, for each test problem,
the new proposed algorithms need almost the same CPU time.

More importantly, we can see from the fifth and sixth columns of the table that
the CPU time of JBDQR is almost 300 to 500 times that of {hyb-CGME} and {hyb-TCGME}
for the test problems shaw and {deriv2} with all noise levels.
And for the test problem heat with the higher noise levels $\varepsilon=10^{-1}$ and $5\times 10^{-2}$,
the CPU time of JBDQR is around 200 to 300 times that of {hyb-CGME} and {hyb-TCGME},
however, with the lower noise level $\varepsilon=10^{-2}$,
it is even 500 to 700 times.
Especially for the severely ill-posed test problem {baart},
the CPU time of JBDQR is around 400 to 600 times that of {hyb-CGME} and {hyb-TCGME}
when the noise levels are $10^{-1},5\times 10^{-2},\ {\rm and } \ 10^{-2}$.
However, for two-dimensional problems, 
We can observe similar phenomena, 
but not as obvious as the one-dimensional problems.
For the test problem {grain} with all noise levels,
the CPU time of JBDQR is almost 10 times that of {hyb-CGME} and {hyb-TCGME}.
For the test problems {mri} and {satellite},
JBDQR takes six to eight times more CPU time than {hyb-CGME} and {hyb-TCGME} do,
and for the test problem {blur440},
the CPU time of JBDQR is around two to three times that of {hyb-CGME} and {hyb-TCGME},
as shown in Table \ref{tab3}.


\begin{figure}

\begin{minipage}{0.48\linewidth}
  \centerline{\includegraphics[width=6.0cm,height=4cm]{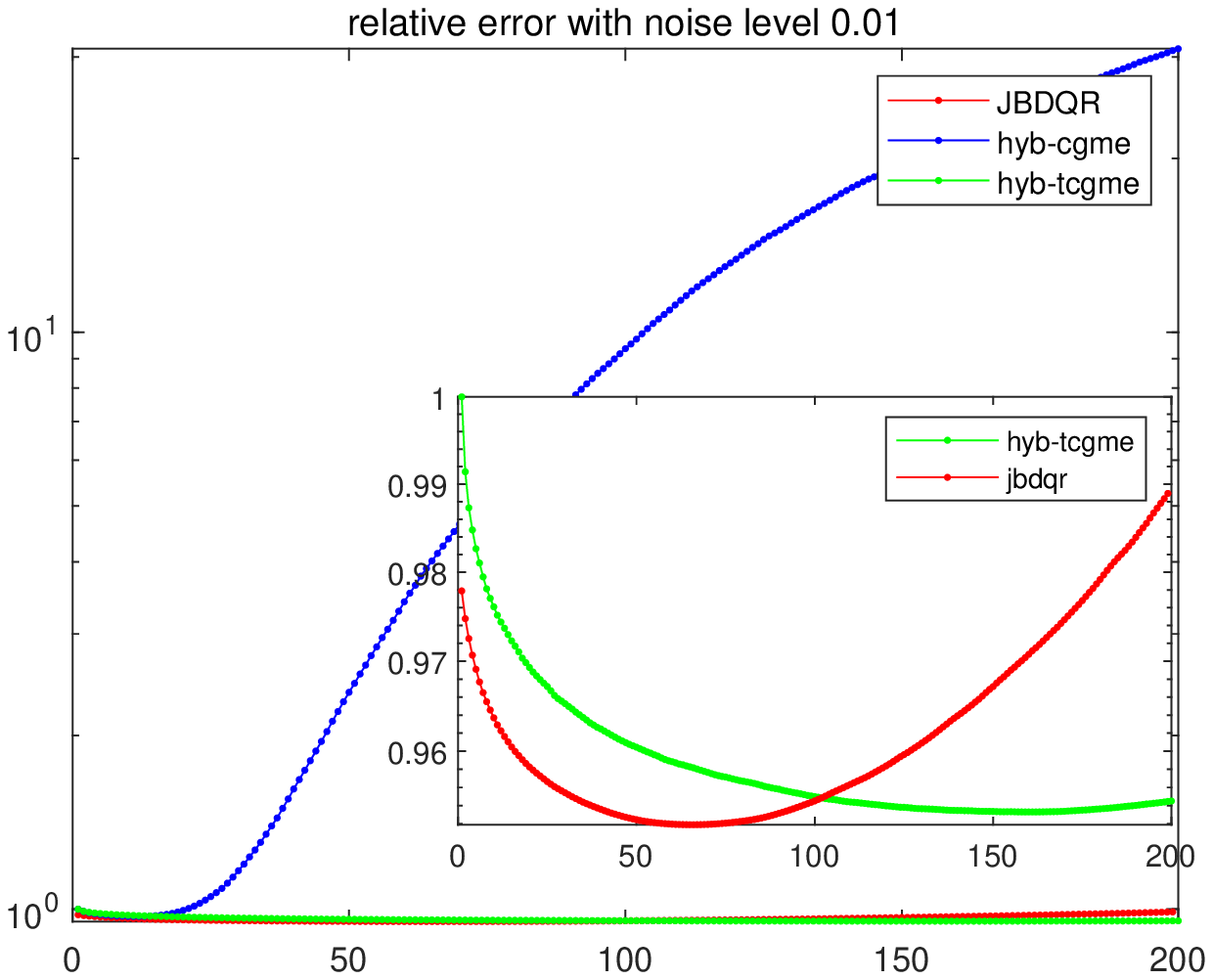}}
  \centerline{(a)}
\end{minipage}
\hfill
\begin{minipage}{0.48\linewidth}
  \centerline{\includegraphics[width=6.0cm,height=4cm]{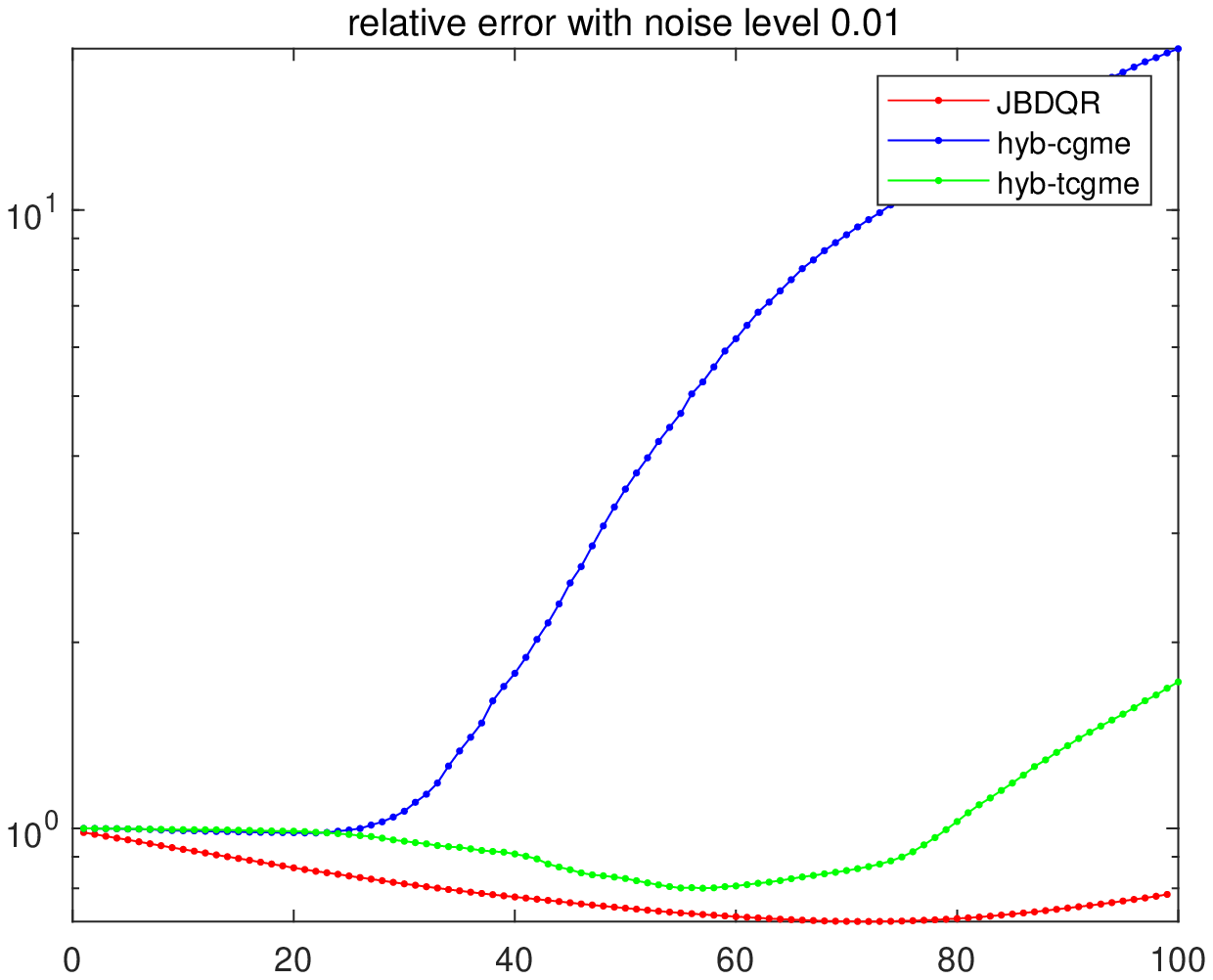}}
  \centerline{(b)}
\end{minipage}
\vfill
\begin{minipage}{0.48\linewidth}
  \centerline{\includegraphics[width=6.0cm,height=4cm]{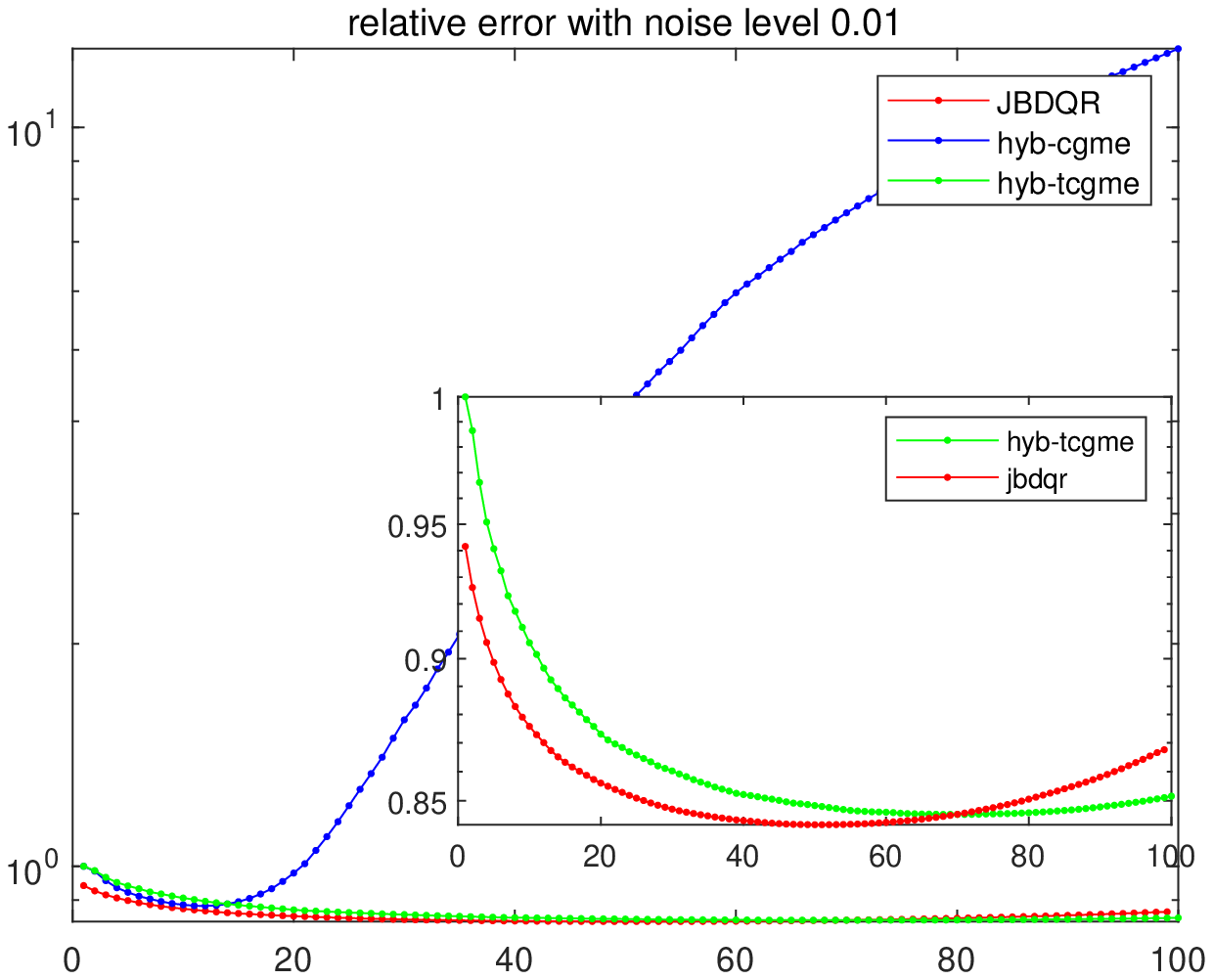}}
  \centerline{(c)}
\end{minipage}
\hfill
\begin{minipage}{0.48\linewidth}
  \centerline{\includegraphics[width=6.0cm,height=4cm]{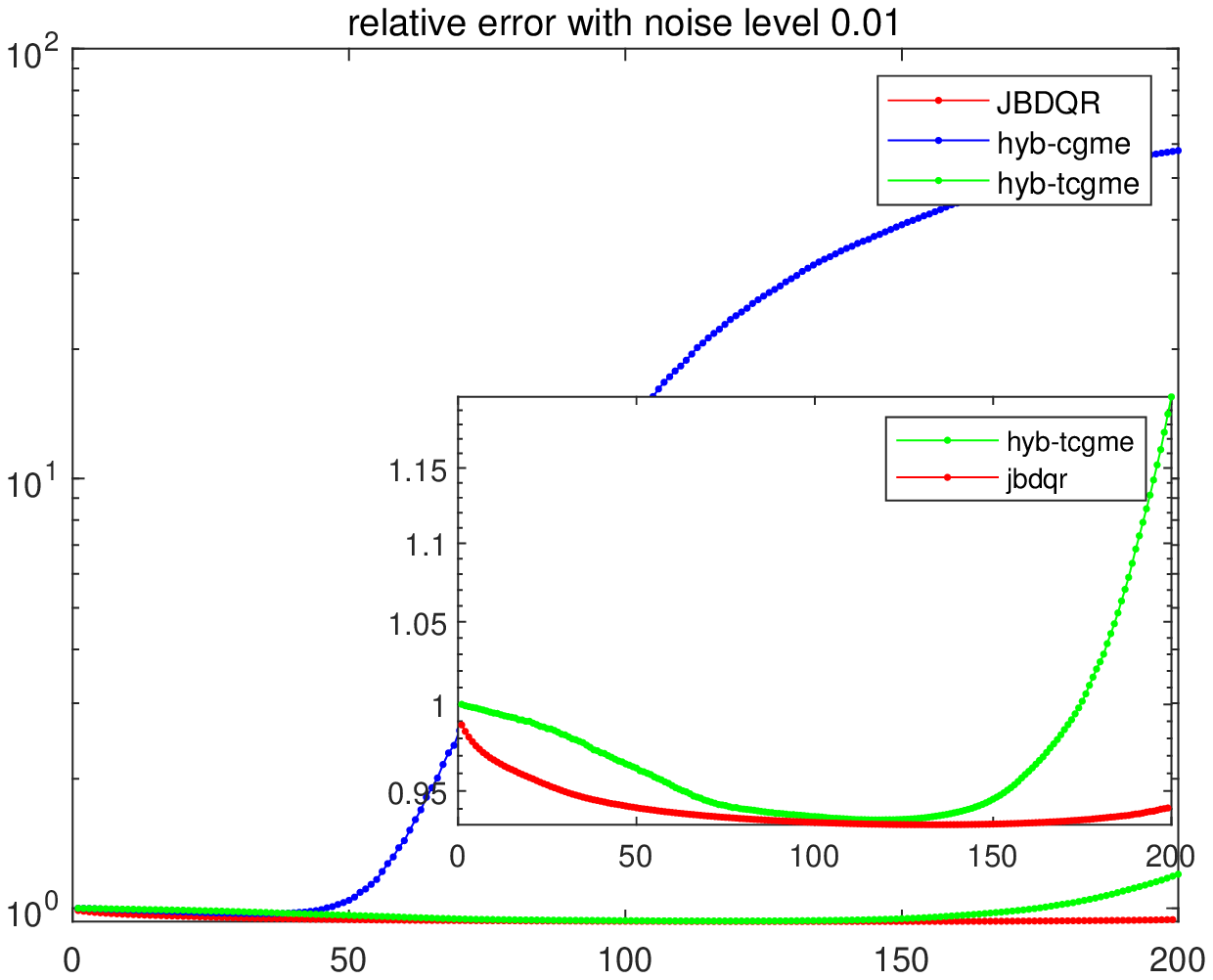}}
  \centerline{(d)}
\end{minipage}
\caption{The relative error of {hyb-CGME}, JBDQR and {hyb-TCGME}
 with $\varepsilon=10^{-2}$: 
(a) {blur440}; (b) {grain}; (c) {mri}; (d) {satellite}.}
\label{fig2}
\end{figure}

\begin{figure}
\begin{minipage}{0.48\linewidth}
  \centerline{\includegraphics[width=6.0cm,height=4cm]{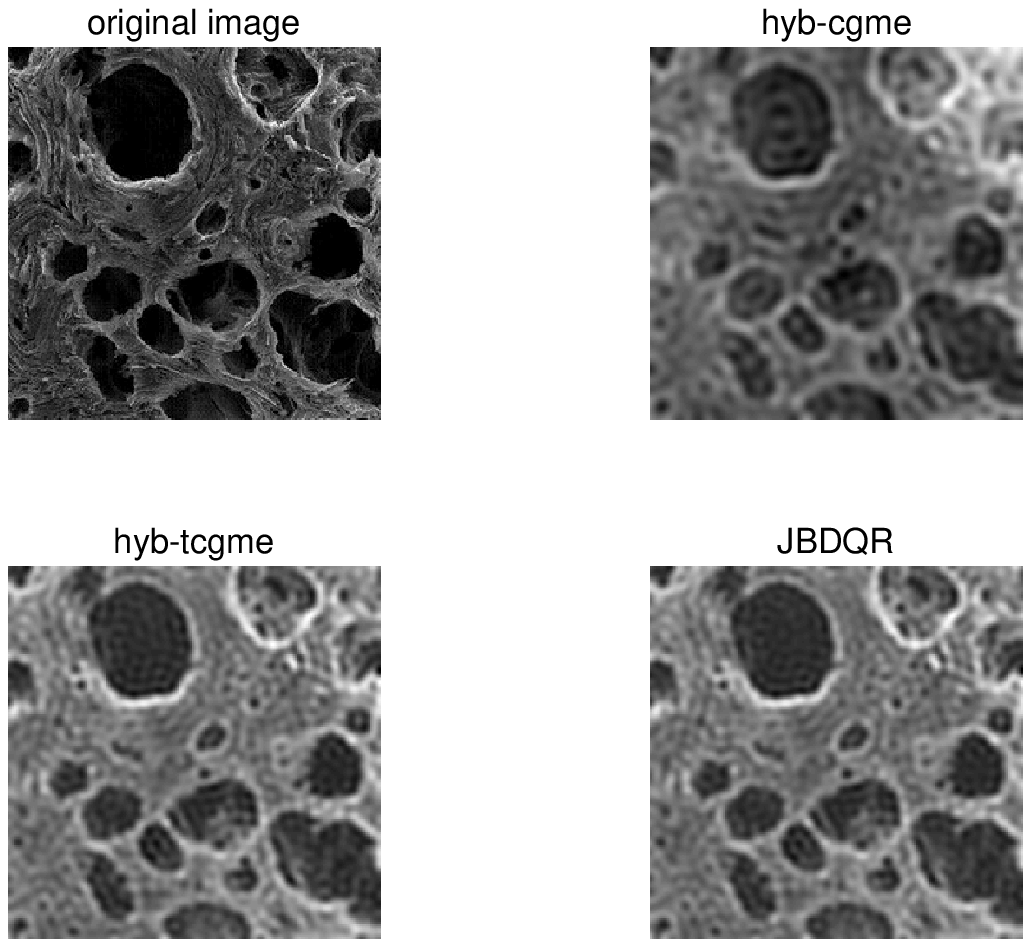}}
  \centerline{(a)}
\end{minipage}
\hfill
\begin{minipage}{0.48\linewidth}
  \centerline{\includegraphics[width=6.0cm,height=4cm]{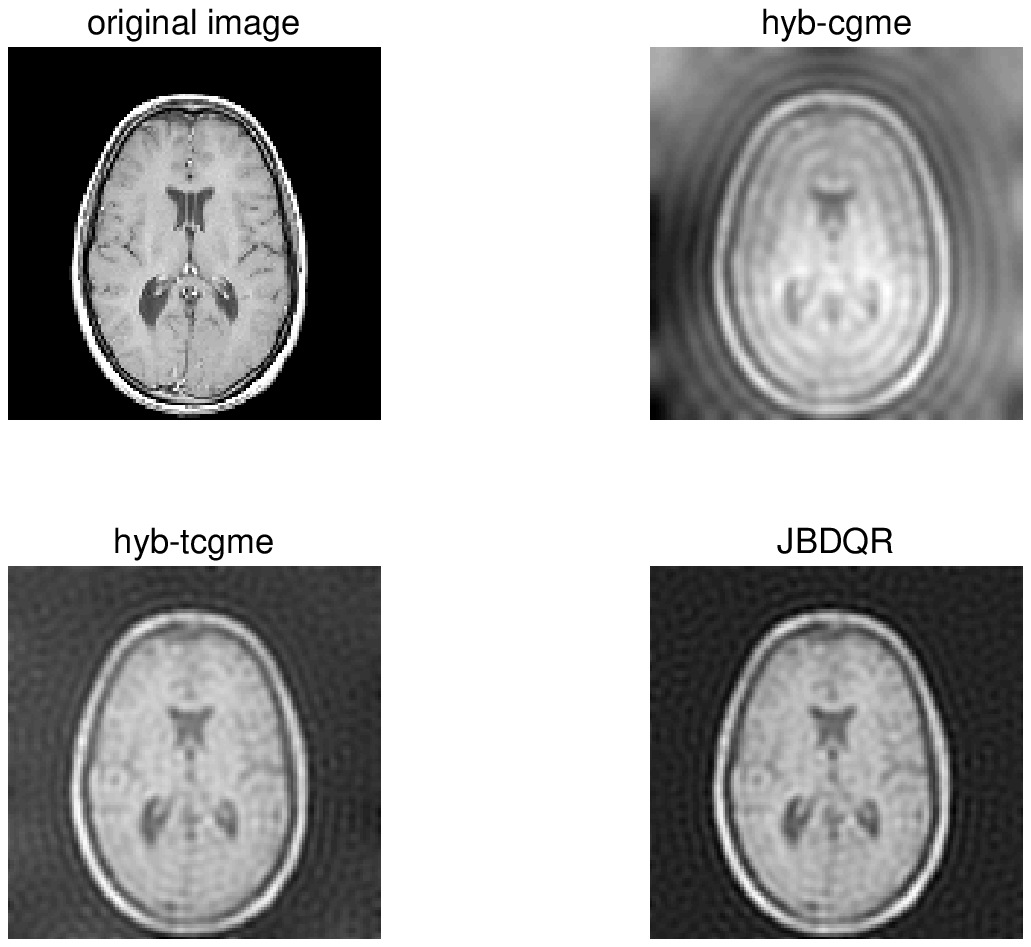}}
  \centerline{(b)}
\end{minipage}
\vfill
\begin{minipage}{0.48\linewidth}
  \centerline{\includegraphics[width=6.0cm,height=4cm]{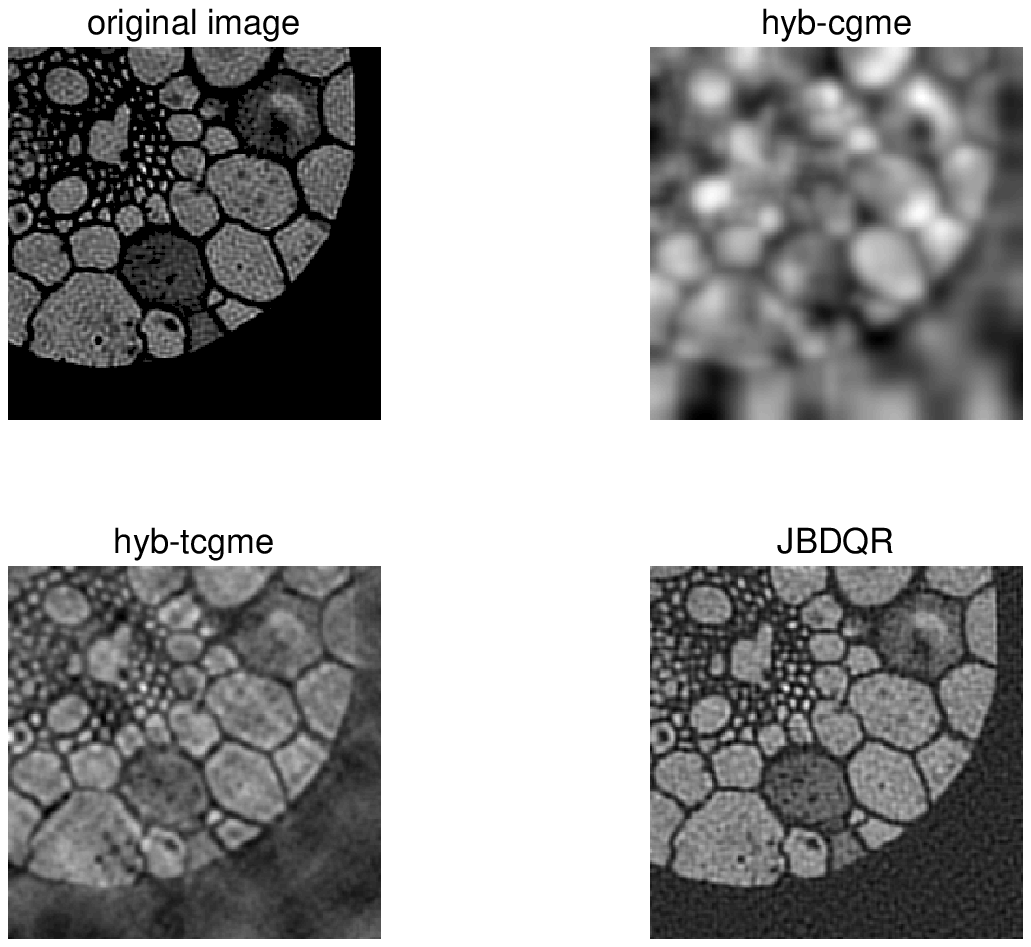}}
  \centerline{(c)}
\end{minipage}
\hfill
\begin{minipage}{0.48\linewidth}
  \centerline{\includegraphics[width=6.0cm,height=4cm]{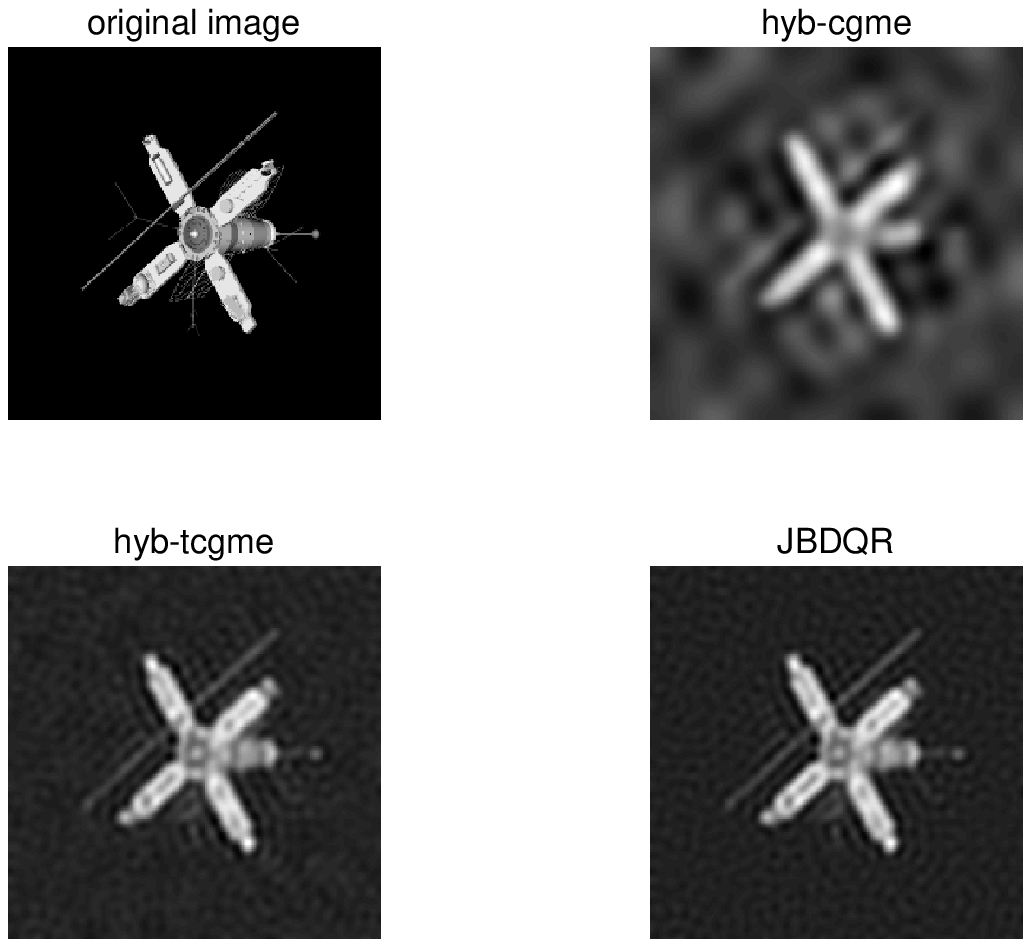}}
  \centerline{(d)}
\end{minipage}
\caption{The exact images and the reconstructed images for the four two dimensional test problems
with $\varepsilon=10^{-2}$ and $L$ defined in (\ref{l3}):
(a) {blur440}; (b) {mri}; (c) {grain}; (d) {satellite}.}
\label{fig3}
\end{figure}

For the sake of length, we only show two-dimensional problems in Figures \ref{fig2} and \ref{fig3}.
In Figure \ref{fig2} we display the convergence processes of {hyb-CGME}, JBDQR and {hyb-TCGME} for $\varepsilon=10^{-2}$.
We can see that the best regularized solutions by {hyb-TCGME} are almost same as
the counterparts by JBDQR except {grain} and
the best regularized solution by {hyb-TCGME} and JBDQR
is better than that by {hyb-CGME} for the four test problems.
Moreover, as the figure shows, for every test problem,
all three algorithms under consideration exhibit semi-convergence \cite{hansen98, hansen10}:
the convergence curves of the three algorithms first decrease with $k$,
then increase.
This means the iterates converge to $x_{true}$
in an initial stage; afterwards the iterates start to diverge from $x_{true}$.
It is worth mentioning that it is proved \cite{jia2020} that the
JBDQR iterates take the form of ﬁltered GSVD
expansions, which shows that JBDQR have the semi-convergence property.
These results indicate that all three algorithms have the semi-convergence property.
We have numerically veriﬁed this assertion, and the details are omitted.


The exact images and the
reconstructed images for the four test problems
with $\varepsilon=10^{-2}$ and $L$ defined by (\ref{l3})
are displayed in Figure \ref{fig3}.
Clearly, the reconstructed
images by {hyb-TCGME} are at least as sharp as those
by JBDQR and the reconstructed
images by {hyb-TCGME} and JBDQR are much more sharp
than the counterparts by {hyb-CGME} .

\section*{Declarations}

%
\begin{itemize}
\item This study was funded by Zhejiang A and F University.
(No. 203402000401).
\item Not applicable.
\item Data availability
\item Materials availability
\item Code availability
\end{itemize}
%



\bibliography{reference}

\end{document}